\documentclass[12pt]{amsart}
\usepackage{amsmath,amsfonts,amsthm,amscd,textcomp,times, amssymb}
 \ifx\pdfoutput\undefined
   \usepackage[dvips]{graphicx}
   \else
   \usepackage[pdftex]{graphicx}
   \fi
  \usepackage{epstopdf}
\newtheorem{theorem}{Theorem}
\newtheorem{conjecture}{Conjecture}
\newtheorem{claim}{Claim}
\newtheorem{proposition}{Proposition}
\newtheorem{lemma}{Lemma}
\newtheorem{definition}{Definition}
\newtheorem{example}{Example}
\newtheorem{corollary}{Corollary}
\newtheorem{remark}{Remark}
\numberwithin{equation}{section}
\numberwithin{theorem}{section}
\numberwithin{proposition}{section}
\numberwithin{lemma}{section}
\numberwithin{claim}{section}
\numberwithin{corollary}{section}
\input{diagrams2010.tex}
\newcommand{\bull}{\ensuremath{{}\bullet{}}}
 \newcommand{\gr}{\ensuremath{\mathbb{G}(N-n,\mathbb{C}^{N+1})}}
\newcommand{\cpn}{\ensuremath{\mathbb{P}^{N}}}
\newcommand{\slnc}{\ensuremath{SL(N+1,\mathbb{C})}}
\newcommand{\dlb}{\ensuremath{\overline{\partial}}}
\newcommand{\dl}{\ensuremath{\partial}}
\newcommand{\ra}{\ensuremath{\longrightarrow}}
\newcommand{\om}{\ensuremath{\omega}}
\newcommand{\vp}{\ensuremath{\varphi}}
\newcommand{\vps}{\ensuremath{\varphi_{\sigma}}}
\newcommand{\oms}{\ensuremath{\omega_{\sigma}}}

\newcommand{\vpt}{\ensuremath{\varphi_{t}}}
\newcommand{\vplt}{\ensuremath{\varphi_{\lambda(t)}}}

\newcommand{\cn}{\ensuremath{\mathbb{C}^{N+1}}}

\newcommand{\xhyp}{\ensuremath{X\times\mathbb{P}^{n-1}}}
\begin{document}
\DeclareGraphicsExtensions{.pdf,.gif,.jpg}

%%%%%%%%%%%%%%%%%%%%%%%%%%%%%%%%%%%%%%%%%%%%%%%%%%%%%%%%%%%%%%%%%%%%%%%%%%%%%%%%%%%%%%%%%%%%%%%%%%%%%%%%%%%%%%%%%%
 \title[Hyperdiscriminant Polytopes ]{Hyperdiscriminant Polytopes, Chow Polytopes, and Mabuchi Energy Asymptotics }
 \author{Sean Timothy Paul}
 \email{stpaul@math.wisc.edu}
 \address{Mathematics Department at the University of Wisconsin, Madison}
 \subjclass[2000]{53C55}
 \keywords{Discriminants, resultants, K-energy maps, projective duality, K\"ahler Einstein metrics, convex polytopes, K-stability  .}
 %%%%%%%%%%%%%%%%%%%%%%%%%%%%%%%%%%%%%%%%%%%%%%%%%%%%%%%%%%%%%%%%%%%%%%%%%%%%%%%%%%%%%%
%\title{ Projective Duality and K-Energy Asymptotics  III}
%\author{Sean Timothy Paul}
%%%%%%%%%%%%%%%%%%%%%%%%%%%%%%%%%%%%%%%%%%%%%%%%%%
\date{September 2, 2010}
 \vspace{-5mm}
\begin{abstract}{Let  $X^n\rightarrow \cpn$ be a smooth, linearly normal algebraic variety. It is shown that the Mabuchi energy of $(X,{\om_{FS}}|_X)$ restricted to the Bergman metrics is completely determined by the $X$-hyperdiscriminant of format $(n-1)$ and the Chow form of $X$. As a corollary it is shown that the Mabuchi energy is bounded from below for all degenerations in $G$ if and only if the hyperdiscriminant polytope dominates the Chow polytope for all maximal algebraic tori $H$ of $G$.}
\end{abstract}
\maketitle
 \tableofcontents \newpage
%%%%%%%%%%%%%%%%%%%%%%%%%%%%%%%%%%%%%%%%%%%%%%%%%%%%%%%%%%
 %%%%%%%%%%%%%%%%%%%%%%%%%INTRO%%%%%%%%%%%%%%%%%%%%%%%%%
\section{Introduction and Statement of Results}
Let $X^n\ra \cpn$ be a smooth complex projective variety of degree $d\geq 2$ embedded by a very ample complete linear system. Fix any Hermitian metric on $\mathbb{C}^{N+1}$ and let $\om_{FS}$ denote the associated Fubini-Study K\"ahler form. We set ${\om:=\om_{FS}}|_X$ . To $\sigma \in G$ ( the automorphism group of $\cpn$ ) we associate the \emph{Bergman potential} $\varphi_{\sigma}\in C^{\infty}(X)$
\begin{align*}
\sigma^*\om_{}=\om_{}+\frac{\sqrt{-1}}{2\pi}\dl\dlb \vps >0 \ .
\end{align*}    
 Let $\nu_{\om}$ denote the Mabuchi energy of $(X,\om)$.  For any $\sigma\in G$ we define 
 \begin{align*}
 \nu_{\om}(\sigma):=\nu_{\om}(\vps) \ .
 \end{align*}
Let $\lambda:\mathbb{C}^*\ra G$ be an algebraic one parameter subgroup of $G$. We shall refer to such maps, and their associated potentials $\varphi_{\lambda(t)}$ , as \emph{degenerations}.   Three basic problems in the field of K\"ahler geometry are the following.\\
 \ \\
 \noindent\textbf{Problem 1.}   \emph{Give a \textbf{complete} description of the behavior of the Mabuchi energy along all degenerations. That is, describe}
    \begin{align*}
  \lim_{|t|\ra 0}  \nu_{\om}(\lambda(t)) \quad  t\in \mathbb{C}^*  \ .
    \end{align*}
  \textbf{Problem 2.} \emph{Provide necessary and sufficient conditions in terms of the geometry of the embedding $X\ra \cpn$ which insure that
 $\nu_{\om}$ is \textbf{bounded below} along all degenerations.}\\
 \ \\
  \textbf{Problem 3.} \emph{Provide necessary and sufficient conditions in terms of the geometry of the embedding which insure that $\nu_{\om}$ is \textbf{proper} along all degenerations.}\\
  
  In this paper we provide complete solutions to all of these problems.  The author's solution is given in terms of the $X$-\emph{resultant} (the Cayley-Chow form of $X$) and the $X$-\emph{hyperdiscriminant} of format $(n-1)$ (the defining polynomial of the variety of tangent hyperplanes to $X\times \mathbb{P}^{n-1}$ in the Segre embedding ). That the  $X$-resultant  appears in the K-energy is not new and is due to Gang Tian (see \cite{kenhyp}) . The author's original contribution is the discovery that the $X$-hyperdiscriminant \emph{also} appears  in the Mabuchi energy of an algebraic manifold . In fact, it is the hyperdiscriminant  that \emph{reflects the presence of the Ricci curvature}. The Chow form does not.  \\
\ \\
\textbf{Theorem A .} 
\emph{ Let $X^n \hookrightarrow \cpn$ be a smooth, linearly normal complex algebraic variety of degree $d \geq 2$ . 
 Let   $R_{X}$  denote the \textbf{X-resultant} (the Cayley-Chow form of $X$) . Let $\Delta_{X\times \mathbb{P}^{n-1}}$ denote the \textbf{X-hyperdiscriminant}  of format $(n-1)$ (the defining polynomial for the dual of $X\times \mathbb{P}^{n-1}$ in the Segre embedding )\footnote{We collect all of the basic definitions in section 2 .}.
 Then the Mabuchi energy restricted to the Bergman metrics is given as follows}
\begin{align}\label{bold}
\begin{split}
&\nu_{\om}(\varphi_{\sigma})= {\deg(R_X)}\log  \frac{{||\sigma\cdot\Delta_{X\times\mathbb{P}^{n-1}}||}^{2}}{{||\Delta_{X\times\mathbb{P}^{n-1}}||}^{2}} - {\deg(\Delta_{X\times\mathbb{P}^{n-1}})} \log  \frac{{||\sigma\cdot R_{X}||}^{2}}{||R_{X}||^2}   \ .
 \end{split}
\end{align}
\begin{remark}\emph{The Mabuchi energy restricted to $G$ is not manifestly, and most likely not, a convex function .}
\end{remark}

It follows from Theorem A that the asymptotic expansion of the Mabuchi energy along any algebraic one parameter subgroup of $H$ (a maximal algebraic torus of $G$ )\footnote{In this paper $G$ always denotes $\slnc$.}  is completely determined by the \emph{Chow polytope} $\mathcal{N}(R_X)$ and the \emph{hyperdiscriminant polytope} $\mathcal{N}(\Delta_{X\times\mathbb{P}^{n-1}})$  (see (\ref{wtpolytope}) ). We remark that these are compact convex lattice polytopes inside $M_{\mathbb{R}}:= M_{\mathbb{Z}}(H)\otimes_{\mathbb{Z}}\mathbb{R}\cong \mathbb{R}^N$, where $M_{\mathbb{Z}}=M_{\mathbb{Z}}(H)$ denotes the rank $N$ lattice of rational characters of $H$. In the statement of Theorem B below $l_{\lambda}$ denotes the integral linear functional on $M_{\mathbb{R}}$ corresponding to the degeneration $\lambda\in N_{\mathbb{Z}}:= M_{\mathbb{Z}}^{\vee}$ (dual lattice) .\\
\ \\
\textbf{Theorem B .} \emph{There is an asymptotic expansion as $|t|\rightarrow 0$}
\begin{align} 
\begin{split}
 &\nu_{\om}(\lambda(t))= F_P(\lambda)\log|t|^2+O(1) \ , \\
\ \\
&F_P(\lambda):={\deg(R_X)}\mbox{{min}}_{ \{ x\in \mathcal{N}(\Delta_{X\times\mathbb{P}^{n-1}})\}}\ l_{\lambda}(x)- 
 {\deg(\Delta_{X\times\mathbb{P}^{n-1}} )} \mbox{ {min}}_{ \{ x\in \mathcal{N}(R_X)\}}\ l_{\lambda}(x) \ .
\end{split}
\end{align}
\emph{In particular, $\nu_{\om}(\lambda(t))$ has a logarithmic singularity as $|t|\rightarrow 0$, and the coefficient of blow up is an integer.}\\
\ \\
Theorem B provides a complete solution to {Problem 1} .\\
\ \\
\textbf{Theorem C .}
\emph{The Mabuchi energy of $(X,{\om_{FS}}|_X)$ is bounded from below along all degenerations in $G$ if and only if for all maximal tori $H$ the hyperdiscriminant polytope dominates the Chow polytope }
\begin{align}\label{bounded}
\deg(\Delta_{X\times\mathbb{P}^{n-1}} )\mathcal{N}(R_X)\subseteq  {\deg(R_X)}\mathcal{N}(\Delta_{X\times\mathbb{P}^{n-1}}) \ .
\end{align}
\ \\
Theorem C provides a complete solution to {Problem 2} .\\
\ \\
\noindent\textbf{Theorem D .}
\emph{The Mabuchi energy of $(X,{\om_{FS}}|_X)$ is {\textbf{proper}} along all degenerations in $G$ if and only if for all
  maximal tori $H$ and all $m>>0 , m\in\mathbb{Z}$ we have}
\begin{align}\label{properthm}
(m- 1)\deg(\Delta_{X\times\mathbb{P}^{n-1}} )\mathcal{N}(R_X)+ \deg(\Delta_{X\times\mathbb{P}^{n-1}} )\deg(R_X)\mathcal{S}_{N}\subseteq m{\deg(R_X)}\mathcal{N}(\Delta_{X\times\mathbb{P}^{n-1}})\ .
 \end{align}
\emph{$\mathcal{S}_{N}$ is the standard $N$-simplex in $\mathbb{R}^N$ and the addition on the left side of (\ref{properthm}) denotes {Minkowski summation} of polyhedra .}\\
\ \\
 Theorem D provides a complete solution to problem 3.
 
 The next result provides a weak form of the {numerical criterion} for the Mabuchi K-energy map. \newpage
\noindent\textbf{Theorem E .}
\emph{Let $H$ be any maximal algebraic torus of $G$. Assume that there is a sequence $\{\tau_i\}\subset H$ such that}
\begin{align*}
\liminf_{i\ra\infty}\nu_{\om}(\varphi_{\tau_i})=-\infty \ .
\end{align*}
\emph{Then there exists a one parameter subgroup $\lambda:\mathbb{C}^*\ra H$ such that}
\begin{align*}
\lim_{|t|\ra 0}\nu_{\om}(\lambda(t))=-\infty \ .
\end{align*}

It seems to be a tacit assumption that the Mabuchi energy is bounded below \emph{generically}. The next result provides a precise quantitative statement to that effect in the context of algebraic degenerations induced from an arbitrary projective embedding.\\
\ \\
\noindent\textbf{Theorem F .}
\emph{Fix a maximal algebraic torus $T$ of $G$. Then there is an explicitly computable algebraic hypersurface $Z=Z({T})\subset G$ such that for all $\sigma\in G\setminus Z$ there is a positive constant $C(\sigma)$ such that}
\begin{align}
\nu_{\om}(\varphi_{\tau})\geq -C(\sigma) \quad \mbox{\emph{for all} $\tau\in \sigma T\sigma^{-1}$} \ .
\end{align}
\ \\
Applications of Theorem A to canonical K\"ahler metrics are as follows, the precise definition of K-(semi)stability is new and appears below (see definitions (\ref{kstability}) and (\ref{fundamental})) .
\begin{corollary}\label{ke}\ \\
  \emph{i) If a polarized manifold $(X,L)$ admits a metric of constant scalar curvature in the class $c_1(L)$ then it is K-semistable with respect to all embeddings  $X\overset{L^m}{\ra}\mathbb{P}^{N_m}$}. \\
\ \\
   \emph{ii)  In particular a Fano manifold  $(X,-K_X)$ admits a K\"ahler Einstein metric only if all pluri-anticanonical models are K-semistable .} \\
 \ \\
 \emph{ iii)  If $(X,-K_X)$ has a discrete symmetry group and admits a K\"ahler Einstein metric then it is K-\textbf{stable}} . 
\end{corollary}
 
 We single out the following special cases.
 
 \begin{corollary}\label{examples}\ \newline\emph{
i)  Any canonically polarized manifold $(X,K_X)$ is K-\textbf{stable} with respect to its pluricanonical embeddings. \\
\ \\
ii) Any polarized Calabi-Yau manifold $(X,L)$ is K-\textbf{stable} with respect to all embeddings  $X\overset{L^m}{\ra}\mathbb{P}^{N_m}$. \\
\ \\
iii)  Any compact homogeneous K\"ahler manifold is K-semistable with respect to its plurianticanonical embeddings }.\\
\end{corollary} 
  
  It should be noted that there is no error term $\Psi$ in (\ref{bold}), unlike the main results of Tian in \cite{kenhyp} (Theorem 4.1 pg. 258), \cite{psc} (formula (8.16) pg. 34), as well Tian and the author (see \cite{ags} pg. 2564 Theorem 3.5). In particular the Mabuchi energy restricted to the Bergman metrics is \emph{not} a  singular, or ``degenerate" norm of the Cayley-Chow form of $X$, but simply the \emph{difference} of two quite honest norms, one involving $R_X$ and the other $\Delta_{X\times\mathbb{P}^{n-1}}$.   Consequently the approach of the author is quite down to earth and focuses on \emph{concrete} (tangents, secants, Gauss maps, etc.) projective geometric constructions with  sub{varieties}  (not schemes) of $\cpn$ very much in the spirit of F.L. Zak \cite{zak} and the seminal paper of Griffiths and Harris \cite{griffharr}.   

Hilbert schemes, test configurations and flat families, limit cycles, generalities concerning $G$-linearized line bundles, Deligne pairings and ``degenerate semi-norms", deformations to the normal cone, and the menagerie of numerical slope stability conditions that are ubiquitous in the literature are all completely unnecessary in the present article.  

\begin{center}\emph{The new perspective in this paper is that the generalized Futaki invariant should not be considered as a number, but rather be interpreted as a \textbf{pair of polytopes} associated to any smooth, linearly normal projective algebraic variety $X\hookrightarrow \cpn$}. \end{center}

The polytopes in question are the \emph{hyperdiscriminant} and \emph{Chow} polytopes of Cayley and Gelfand, Kapranov and Zelevinsky (see \cite{newtpoly}, \cite{gkz} and \cite{ksz}) . The {test configurations} in the literature on K-stability are {linear functionals} on these polytopes. The difference of the minima of these functionals is what controls the K-energy map for any smooth algebraic variety.  From the author's new point of view \emph{degenerating the variety is not necessary}. The problem is to {understand the relative positions of these polytopes}. This does not require full knowledge of the $X$-resultant and $X$-discriminant, only their  {{supports}} are relevant.   We should point out that {our expression (\ref{bold}) for the K-energy map is given for \emph{all} the Bergman metrics, not merely the diagonal ones}.     

  Theorem B provides a new definition of the {generalized Futaki invariant} that (i) does not require smoothness (or normality) of the limit cycle and (ii) completely captures the behavior of the Mabuchi energy along the degeneration.  In the case of a smooth limit cycle, our definition agrees with the original definition of Ding and Tian . The generalized Futaki invariant proposed by Donaldson in \cite{skdtoric} satisfies (i) but only satisfies (ii) in the special case of reduced limit cycle.  Consequently, our work gives a new understanding of {K-stability} . 

Although the table of contents reveals the organization of the paper, a few remarks may be helpful. The proof Theorem A is contained in sections 4 and 5. The remaining results follow from the discussion in section 2. The reader is advised to first consult definitions \ref{kstability}, \ref{fundamental} and glance through section 3  and then proceed directly to sections 4 and 5.  

\subsection{Notations and Preliminaries}%%%%%%%%%%%%%%%%%%%%%%%%%%%%%%%%%%%%%%%%%%%%%%%%%%%%%%%%%
Let $(X,\ \om)$ be a K\"ahler manifold. We always set $\mu$ to be the average of the scalar curvature of $\om$ and $V$ to be the volume
\begin{align*}
\mu:=\frac{1}{V}\int_X\mbox{Scal}(\om)\om^n \ , \
V:= \int_X\om^n \ .
\end{align*}

The space of K\"ahler potentials will be denoted by $\mathcal{H}_\om$
\begin{align*}
\mathcal{H}_\om:=\{\vp\in C^{\infty}(X)|\om_\vp:=\om+\frac{\sqrt{-1}}{2\pi}\dl\dlb\vp>0 \} \ .
\end{align*}

The \emph{Mabuchi K-energy}, denoted by $\nu_\om$,  is a map $\nu_\om:\mathcal{H}_\om\ra \mathbb{R}$ and is given by the following expression \footnote{{We warn the reader that our definition of $\nu_{\om}$ differs from the usual one by a factor of $V^2(n+1)$}. }
\begin{align}\label{mabenergy}
 \qquad \nu_{\omega}(\varphi):= -(n+1){V}\int_{0}^{1}\int_{X}\dot{\varphi_{t}}(\mbox{Scal}(\varphi_{t})-\mu)\omega_{t}^{n}dt.
\end{align}
Above, $\varphi_{t}$ is a smooth path in $\mathcal{H}_\om$ joining $0$ with $\varphi$. It is well known that the K-energy does not depend on the path chosen (see \cite{mabuchi}) . $\vp$ is a critical point of the Mabuchi energy  if and only if $\mbox{Scal}_\om(\vp)\equiv \mu$ . 
Suppose that $\om$ satisfies $\mbox{Ric}(\om)= \frac{\mu}{n}\om+\frac{\sqrt{-1}}{2\pi}\dl\dlb h_\om$ .   In this case there is the following well known direct formula for the K-energy map. 

\begin{align} \label{directformula}
\begin{split}
& \nu_{\omega}(\varphi)=\kappa\Big\{\frac{1}{V}\int_{X}\mbox{log}\left(\frac{{\omega_{\varphi}}^{n}}{\omega^{n}}\right){{\omega_{\varphi}}^{n}} - \frac{\mu}{n}(I_{\omega}(\varphi)-J_{\omega}(\varphi)) - \frac{1}{V}\int_Xh_\om(\om^n_\vp - \om^n)\Big\} \\
&\kappa:=V^2(n+1) \ . \\
\ \\
& J_{\omega}(\varphi):= \frac{1}{V}\int_{X}\sum_{i=0}^{n-1}\frac{\sqrt{-1}}{2\pi}\frac{i+1}{n+1}\dl\varphi \wedge \dlb
\varphi\wedge \omega^{i}\wedge {\omega_{\varphi} }^{n-i-1}\\
\ \\
&I_{\omega}(\varphi):= \frac{1}{V}\int_{X}\varphi(\omega^{n}-{\omega_{\varphi}}^{n}) \ .
\end{split}
\end{align}

\begin{definition}\label{proper}
\emph{The Mabuchi energy is \emph{\textbf{proper}} provided there exists constants $A>0$ and $B>0$ such that for all $\varphi\in \mathcal{H}_\om$ we have }
\begin{align}
\begin{split}
&\nu_{\om}(\vp)\geq A J_{\omega}(\varphi)-B  \ .
 \end{split}
\end{align}
\end{definition}

The notion of properness is due to Tian (see \cite{psc}).  Assume that $\om$ represents a multiple of the canonical class.  Jensen's inequality applied to (\ref{directformula}) together with the simple inequality
\begin{align}
\frac{1}{n+1}I_{\om}(\vp)\leq J_{\omega}(\varphi)\leq \frac{n}{n+1}I_{\om}(\vp)
\end{align}
shows that $\nu_{\omega}$ is proper whenever $\mu\leq 0$ .

Two basic facts that we require in this paper are the following.
 
\begin{theorem}[Bando, Mabuchi \cite{bandmab}]\label{bandmab} Let $(X,\ \om)$ be a compact Fano manifold with $[\om]=c_1(X)$. If there exits a K\"ahler Einstein metric in the class $[\om]$ then $\nu_{\omega}$ is bounded from below.
\end{theorem}

 \begin{theorem}[Tian \cite{psc}]\label{tian97} Let $(X,\ \om)$ be a compact Fano manifold with $[\om]=c_1(X)$. Assume that $\eta(X)=\{0\}$. Then there exits a K\"ahler Einstein metric in the class $[\om]$ if and only if $\nu_{\omega}$ is proper.
  \end{theorem}
     
%%%%%%%%%%%%%%%%%%%%%%%%%%%%%%%%%%%%%%%%%%%%%%%%%%%%%%%%%%%%%%%%%%%%%%%%%%%%%%%%%%%%%%%%%%%%%%%%%%%%%%%%%%%%%%%%%%%%%%%%%%%%%%%%%%%%%%%%%%%%%%%%%%%%%%%%%%%%%%%%%%%%%%%%%%%%%%
\section{K-stability and representations of the special linear group}
Let $G$ be one of the classical subgroups of $GL(N+1,\mathbb{C})$. For the most part we shall consider the case
\begin{align*}
G=\slnc \ .
\end{align*}
 Let $(V,\rho)$ be a finite dimensional complex rational representation of $G$. Recall that $E$ is \emph{rational} provided that for all $\alpha\in V^{\vee}$ (dual space) and $v\in V\setminus \{0\}$ the \emph{matrix coefficient}
\begin{align*}
\varphi_{\alpha , v}:G\ra \mathbb{C} \quad \varphi_{\alpha , v}(\sigma):=\alpha(\rho(\sigma)\cdot v) 
\end{align*}
is a \emph{regular function} on $G$, that is
\begin{align*}
\varphi_{\alpha , v}\in \mathbb{C}[G]:= \mbox{affine coordinate ring of $G$} .
\end{align*}
To begin, let $H$ denote any maximal algebraic torus of $G$.  $M_{\mathbb{Z}}=M_{\mathbb{Z}}(H)$ denotes the \emph{character lattice} of $H$
\begin{align*}
M_{\mathbb{Z}}:= Hom_{\mathbb{Z}}(H,\mathbb{C}^*) \ . 
\end{align*}
$M_{\mathbb{Z}}$ consists of algebraic homorphisms $\chi:H\ra \mathbb{C}^*$ . If we fix an isomorphism 
\begin{align}
M_{\mathbb{Z}}\cong \mathbb{Z}^N 
\end{align}
then we may express each such $\chi$ as a \emph{Laurent monomial}
\begin{align*}
\chi  (t_{1},t_{2},\dots,t_{N})=t_{1}^{m_{1}}t_{2}^{m_{2}}\dots t_{N}^{m_{N}} \ , \quad m_{i} \in \mathbb{Z}\ .
\end{align*}
Therefore we make the identification
\begin{align*}
\chi = (m_1,m_2,\dots,m_N)\in \mathbb{Z}^N \ . 
\end{align*}

We denote the dual lattice by $N_{\mathbb{Z}}$. It is well known that $N_{\mathbb{Z}}$ consists of the algebraic one parameter subgroups $\lambda$ of $H$. These are algebraic homomorphims $\lambda:\mathbb{C}^*\ra H$.  The duality 
is given by 
\begin{align}
<\cdot \ , \ \cdot>:N_{\mathbb{Z}}\times M_{\mathbb{Z}}\ra \mathbb{Z} \ , \ \chi(\lambda(t))=t^{<\lambda , \chi>} \ .
\end{align}

Since we have fixed some isomorphism of $H$ with the standard torus in $G$ we have    
\begin{align*}
\lambda(t)=\begin{pmatrix}t^{n_{1}}&\dots&\dots& 0\\
                             0&t^{n_{2}}&\dots& 0\\
                             0&\dots&\dots& t^{n_{N}}
                             \end{pmatrix}\ .\end{align*}
 
In this case the pairing is given concretely as follows
\begin{align*}
<\lambda,\chi>=m_1n_1+m_2n_2+\dots +m_Nn_N \ .
\end{align*}
 We introduce the corresponding real vector spaces by extending scalars 
\begin{align*} 
& M_{\mathbb{R}}:= M_{\mathbb{Z}}\otimes_{\mathbb{Z}}\mathbb{R}\cong \mathbb{R}^N\\
\ \\
&N_{\mathbb{R}}:=N_{\mathbb{Z}}\otimes_{\mathbb{Z}}\mathbb{R}=M_{\mathbb{R}}^{\vee} \ .
\end{align*}

The image of $\lambda$ in $N_{\mathbb{R}}$ is denoted by $l_{\lambda}$. Then $l_{\lambda}$ is an \emph{integral} linear functional on $M_{\mathbb{R}}$ .
Since $V$ is rational it decomposes under the action of $H$ into \emph{weight spaces}
\begin{align}
\begin{split}
&V=\bigoplus_{\chi\in \mbox{supp}(V)}V_{\chi}  \\
\ \\
& V_{\chi}:=\{v\in V\ |\ h\cdot v=\chi(h) v \ , \ h\in H\}
\end{split}
\end{align}
where we have defined the \emph{support} of $V$ by
\begin{align*}
supp(V):= \{\chi\in M_{\mathbb{Z}}\ | \ V_{\chi}\neq 0\} \ .
\end{align*}
Given $v\in V\setminus \{0\}$  the projection of $v$ into $V_{\chi}$ is denoted by $v_{\chi}$. The support of any (nonzero) vector $v$ is then defined by
\begin{align*}
supp(v):= \{\chi\in M_{\mathbb{Z}}\ | \ v_{\chi}\neq 0\} \ .
\end{align*}
\begin{definition}   Let $H$ be any maximal torus in $G$. Let $v\in V\setminus\{0\}$ . The \textbf{weight polytope} of $v$ is the compact convex integral polytope $\mathcal{N}(v)$ given by
\begin{align}\label{wtpolytope}
\mathcal{N}(v):=\mbox{\emph{convex hull of the lattice points}}\{\chi\in supp(v)\}\subset M_{\mathbb{R}} \ .
\end{align}
\end{definition}
 In the same vein we define the weight polytope of the module itself by
 \begin{align}
 \mathcal{N}(V):=\mbox{\emph{convex hull of the lattice points}}\{\chi\in supp(V)\}\subset M_{\mathbb{R}} \ .
 \end{align}
 Obviously $\mathcal{N}(v)\subseteq \mathcal{N}(V)$ for  any $v\in V\setminus \{0\}$ and $H\leq G$, when \emph{equality} holds we say that $v$  is \emph{\textbf{generic}} with respect to $H$.

 Let $\mathbb{C}^{N+1}$ denote the standard representation of $G$ and let $H$ be a maximal algebraic torus. The \emph{\textbf{standard simplex}} denoted by $\mathcal{S}_N$ is defined to be the weight polytope of any $H$ generic vector $u\in\mathbb{C}^{N+1}\setminus\{0\}$
\begin{align*}
\mathcal{S}_N:= \mathcal{N}(u)\subset M_{\mathbb{R}}  \ .
\end{align*}
This is an $N$-dimensional polytope containing the origin in its interior.  Next fix any $H\leq G$, we define the \emph{\textbf{degree}} $q(V)$ of the representation as follows
\begin{align}
q(V):=\min\Big\{k\in \mathbb{Z}_+ \ |\ \mathcal{N}(V)\subseteq k\mathcal{S}_N\ \Big\} \ .
\end{align}
 
Now we are prepared to introduce our fundamental definition. Below $G=\slnc$ and $H\leq G$ is a maximal algebraic torus .
\begin{definition}\label{kstability}\emph{ Let $V$ and $W$ be finite dimensional complex rational representations of $G$. Let $v\in V\setminus \{0\}$ and $w\in W\setminus \{0\}$ .
\begin{enumerate}
\item The {pair} $(v,w)$ is \textbf{\emph{K-semistable}} with respect to $H$ if and only if $\mathcal{N}(v)\subseteq \mathcal{N}(w)$ .\\
\ \\
\item $(v,w)$ is K-semistable with respect to $G$ if and only if it is K-semistable for all maximal tori $H$ in $G$.   \\
\ \\
\item $(v,w)$ is \textbf{\emph{K-stable}} with respect to $H$ if and only if there exists $m_0\in\mathbb{N}$ such that
\begin{align*}
(v^{\otimes (m-1)}\otimes u^{q(V)}\ , \ w^{\otimes m})
\end{align*}
is K-semistable with respect to $H$ for all $m\geq m_0$ and $u$ is any $H$-generic vector in the standard representation of $G$.
\ \\
\item $(v,w)$ is K-stable with respect to $G$ if and only if it is K-stable for all maximal tori $H$ in $G$.  
\end{enumerate}}
\end{definition}

 When $V$ is \emph{irreducible} it is well known that $V$ is located in a unique tensor power of the standard representation
 \begin{align}
 V\subset (\cn)^{\otimes p}\qquad p\in \mathbb{Z}_+ \ .
 \end{align}
 In this case we have $p=q(V)$ .
 
 That $q$ depends only on $(V,\rho)$ and not on $H$ in the general case follows from
\begin{proposition}\label{adjointrep}
Fix a maximal torus $H$, let $v\in V\setminus\{0\}$, and let $\sigma\in G$. Then we have the relation
\begin{align}
 Ad(\sigma)\Big(\mathcal{N}_{H}(\sigma\cdot v)\Big)=\mathcal{N}_{\sigma^{-1}H\sigma} (v)
\end{align}
Where $Ad(\sigma)$ denotes the linear extension of the induced equivalence of $\mathbb{Z}$ modules
\begin{align}
\begin{split}
&Ad(\sigma):M_{\mathbb{Z}}(H)\overset{\cong}{\ra}   M_{\mathbb{Z}}(\sigma^{-1}H\sigma) \\
\ \\
&Ad(\sigma)(\chi)(\tau):= \chi(\sigma \tau \sigma^{-1}) \ \mbox{for all $\tau\in\sigma^{-1}H\sigma$ }\ .
\end{split}
\end{align}
\end{proposition}
We have formulated K-stability in terms of arbitrary finite dimensional $G$-modules $V$ and $W$. In our main applications the modules are not only both \emph{irreducible} but satisfy further conditions which we will now consider.

To begin, let $\lambda_{\bull}$ be a partition consisting of $N$ parts 
\begin{align}
\lambda_{\bull}=(\lambda_1\geq \lambda_2\geq \dots\geq \lambda_N\geq \lambda_{N+1}=0)\ .
\end{align}
We let $\mathbb{S}_{\lambda_{\bull}}(\cn)$ denote the corresponding irreducible representation of $G$  with highest weight $\lambda_{\bull}$ (with respect to a maximal algebraic torus $H$) . Let $W_G$ denote the Weyl group of $G$ with respect to $H$, then the weight polytope of the module is given by
\begin{align}
\mathcal{N}(\lambda_{\bull})=\mbox{convexhull} \ \big\{W_G\cdot {\lambda_{\bull}}\big\} \ .
\end{align}  
 Where $W_G\cdot {\lambda_{\bull}}$ denotes the orbit of the highest weight under the action of the Weyl group.
 Consider two irreducible $G$ modules $\mathbb{S}_{\lambda_{\bull}}(\cn)$ and $\mathbb{S}_{\mu_{\bull}}(\cn)$ satisfying the following two conditions (we shall say that the partitions are \emph{admissible}) 
\begin{align}
\begin{split}
&i) \ | \lambda_{\bull}|= | \mu_{\bull}|\\
\ \\
&ii) \ \lambda_{\bull} \trianglelefteq \mu_{\bull}  \ .
\end{split}
\end{align}
Where $| \lambda_{\bull}|:= \sum_{j=0}^N\lambda_j$  and $\trianglelefteq$ denotes \emph{dominance order} : 
\begin{align}
\begin{split}
&   \lambda_{\bull} \trianglelefteq \mu_{\bull} \ \mbox{if and only if for all $1\leq i\leq N$ we have }\ \sum_{k=1}^i\lambda_k\leq \sum_{k=1}^i\mu_k \ .
\end{split}
\end{align}
The following proposition seems to be well known 
\begin{proposition} Let $\mathbb{S}_{\lambda_{\bull}}(\cn)$ and $\mathbb{S}_{\mu_{\bull}}(\cn)$ be two irreducible $G$-modules. Assume that $| \lambda_{\bull}|= | \mu_{\bull}|$ then
 \begin{align}
\lambda_{\bull} \trianglelefteq \mu_{\bull} \ \mbox{if and only if}\ \mathcal{N}(\lambda_{\bull})\subseteq \mathcal{N}(\mu_{\bull}) \ .
\end{align}
\end{proposition}
Fix $n\in\mathbb{Z}_+$ and choose $d\in\mathbb{Z}_+$ satisfying $d\equiv 0\ \mbox{mod} \ n(n+1) $ .  Our main application of K-stability involves the following specific highest weights :
\begin{align}\label{highweight}
\begin{split}
&\lambda_{\bull}=\frac{1}{n+1}\big(\overbrace{ {d} , {d} ,\dots, {d}}^{n+1},\overbrace{0,\dots,0}^{N-n}\big) \ \mbox{and}\ \mu_{\bull}=\frac{1}{n}\big( \overbrace{ {d} , {d} ,\dots, {d}}^n ,\overbrace{ 0,\dots,0}^{N+1-n}\big)
\end{split}
\end{align}
Then it is well known that (see \cite{tableaux}) the corresponding irreducible modules are given by:
\begin{align}
\begin{split}
&\mathbb{S}_{\lambda_{\bull}}(\cn)\cong H^0(\mathbb{G}(N-n-1,\cpn),\mathcal{O}(\frac{d}{n+1}))\\
\ \\
&\mathbb{S}_{\mu_{\bull}}(\cn)\cong H^0(\mathbb{G}(N-n,\cpn),\mathcal{O}(\frac{d}{n}))\\
\end{split}
\end{align}

Obviously $|\lambda_{\bull}|=|\mu_{\bull}|$ and $\lambda_{\bull} \trianglelefteq \mu_{\bull}$. In this case we may verify the polytope inclusion directly. To begin let $0<k<l$; define
\begin{align}
\mathcal{A}_{k, l}:=\big\{e_{i_1}+e_{i_2}+\dots +e_{i_k}\ | \ 1\leq i_1<i_2<\dots <i_k\leq l\big\} \ .
\end{align}
Then the \emph{hypersimplex} $\Delta(k,l)$ of type $(k,l)$ is given by
\begin{align}
\Delta(k,l):=\mbox{convex hull}\ \mathcal{A}_{k, l} =\mbox{convex hull}\ \mathsf{S}_l\cdot(e_1+e_2+\dots +e_k)\ .
\end{align}
$\mathsf{S}_l\cong W_{SL(l,\mathbb{C})}$ denotes the symmetric group.
It is clear that when the weights are given by (\ref{highweight}) we have
\begin{align}
\begin{split}
&\mathcal{N}(\lambda_{\bull})=\frac{d}{n+1}\Delta(n+1,N+1) \\
\ \\
&\mathcal{N}(\mu_{\bull})=\frac{d}{n}\Delta(n,N+1) \ .
\end{split}
\end{align}
The inclusion
\begin{align}
n\Delta(n+1, N+1)\subseteq (n+1)\Delta(n,N+1) 
\end{align}
follows at once from the equality
\begin{align}
n(e_{i_1}+e_{i_2}+\dots +e_{i_{n+1}})=\frac{1}{(n+1)}\Big(\sum_{J\subset \{i_1,i_2,\dots,i_{n+1}\}\ |J|=n}(n+1)\sum_{j\in J}e_j\Big) \ .
\end{align}

Let $(\lambda_{\bull}, \mu_{\bull})$ be admissible, then the  K-semistability of a pair $(v,w)$ holds for \emph{generic} maximal algebraic tori. One sees this as follows.  Fix $v\in \mathbb{S}_{\lambda_{\bull}}(\cn)\setminus \{0\}$, where we have fixed a maximal algebraic torus $T$ .  Then we define a non-trivial {polynomial} $Q$ on $G$ as follows
\begin{align*}
Q_{\lambda_{\bull}; \ v}(\sigma):=\prod_{s\in W_G}<\sigma\cdot v,e_{s\cdot \lambda_{\bull}}> \ .
\end{align*}
$< \ , \ >$ denotes any inner product rendering the weight space decomposition of $\mathbb{S}_{\lambda_{\bull}}(\cn)$ under $T$ as an orthogonal decomposition. Once more, $W_G$ denotes the \emph{Weyl group} of $G$ with respect to $T$, and $e_{s\cdot \lambda_{\bull}}$ denotes the (unique up to scale) weight vector corresponding to the image of the highest weight under the action of $s\in W_G$. In particular, $Q_{\lambda_{\bull}; v}$ is given only up to scale. Now we define
\begin{align}
\begin{split}
& Z_{\lambda_{\bull}; \ v}:=\{\sigma\in G\ | \ Q_{\lambda_{\bull}; v}(\sigma)=0\} \\
\ \\
& U_{{\lambda_{\bull}; \ v}}:=G\setminus Z_{\lambda_{\bull}; v} \ .
\end{split}
\end{align}
Then we have the following
\begin{proposition}
 {For all $\sigma\in U_{{\lambda_{\bull}; \ v}}$ we have the equality of polytopes 
\begin{align}
\mathcal{N}(v)=\mathcal{N}(\lambda_{\bull}) \ 
\end{align}
relative to the torus $H=\sigma T \sigma^{-1}$ . }
\end{proposition}
Given two irreducible $G$ modules with admissible highest weights $\lambda_{\bull}$ and $\mu_{\bull}$  respectively  we have the following corollary.
\begin{corollary} \label{generic}
 {Fix a maximal algebraic torus $T$. Let $v\in \mathbb{S}_{\lambda_{\bull}}(\cn)\setminus \{0\}$ and $w\in \mathbb{S}_{\mu_{\bull}}(\cn)\setminus \{0\}$, then for all $\sigma \in U_{{\mu_{\bull}; \ w}}$ the pair $(v,w)$ is K-semistable with respect to $\sigma T \sigma^{-1}$.}
\end{corollary}

Given two regular $G$ modules $V$ and $W$ observe that for any $H\leq G$ we have
 \begin{align}
\begin{split}
& \mathcal{N}(v\otimes w)=\mathcal{N}(v)+\mathcal{N}(w) \quad v\otimes w\in V\otimes W \ .
  \end{split}
 \end{align}
\begin{example} [\emph{Relation to Hilbert-Mumford Stability}] 
\emph{The reader may easily verify the following proposition which demonstrates, among other things, that Hilbert-Mumford stability is a \emph{special case} of K-stability. In particular it provides many examples of K-semistable pairs.
\begin{proposition} Let $d\in \mathbb{Z}$, $d\geq 2$. Let $V$ be a rational representation of $G$, $v\in V\setminus \{0\}$. 
\begin{enumerate}
\item $(v,v^{\otimes d})$ is K-semistable if and only if $v$ is Hilbert-Mumford stable in the ordinary sense, that is, $0\notin \overline{G\cdot v}$. \\
\ \\
\item If $(v,v^{\otimes d})$ K-stable if and only if $v$ is (strictly) Hilbert Mumford stable in the ordinary sense, that is $\overline{G\cdot v}={G\cdot v}$ and $G_v$ is finite . \\
\ \\
\item Assume that $(v ,w)$ is K-semistable. If $v $ is Hilbert-Mumford stable then so is $w$.
\end{enumerate}
\end{proposition}}
\end{example}
 \begin{example}[\emph{Classical discriminant and Resultants}]\emph{Consider two polynomials $P$ and $Q$ in one variable of degrees $m$ and $n$ respectively
\begin{align*}
&P(z)=a_mz^m+a_{m-1}z^{m-1}+\dots + a_1z+a_0 \\
\ \\
&Q(z)=b_nz^n+b_{n-1}z^{n-1}+\dots + b_1z+b_0 \ .
\end{align*}
Recall that the classical \emph{resultant} of $P$ and $Q$ is the (quasi)homogeneous polynomial  of the coefficients $(a_0,\dots,a_m;b_0,\dots,b_n)$ defined by
\begin{align*}
R_{m,n}(P,Q)= R_{m,n}(a_0,\dots,a_m;b_0,\dots,b_n):=b_n^m\prod_{\beta_i\in \mbox{zer}(Q)}P(\beta_i)=(-1)^{mn}R_{n,m}(Q,P) \ .
\end{align*}
When $m=n=d\geq 2$ we denote the resultant by $R_d$. Then
\begin{align*}
R_d\in \mathbb{C}_{2d}[M_{2\times (d+1)}] \ .
\end{align*}
$G=SL(d+1,\mathbb{C})$ acts on $R_d$ by the rule
\begin{align*}
\sigma\cdot R_d (A):=R_d(A\cdot \sigma) \quad \sigma\in G \ , \ A\in M_{2\times (d+1)} \ .
\end{align*}
The \emph{discriminant} , $\Delta_d$ , of a polynomial $P$ of degree $d$ is defined by
\begin{align*}
&\Delta_d(a_0,\dots,a_d):=R_{d,d-1}(P,\frac{\dl P}{\dl z}) \\
\ \\
&\Delta_d\in\mathbb{C}_{2d-2}[M_{1\times(d+1)}] \ .
\end{align*}
The action of $G$ is given by
\begin{align*}
\sigma \cdot \Delta_d(a)=\Delta_d(a\cdot \sigma) \ .
\end{align*}
It follows from beautiful work of Gelfand, Kapranov and  Zelevinsky (\cite{newtpoly}) that the pair $(R_d^{\deg(\Delta_d)}, \Delta_d^{\deg(R_d)})$ is K-semistable with respect to the \emph{standard} torus, i.e. the torus corresponding to the $dth$ Veronese embedding of $\mathbb{P}^1$ .
\begin{claim}\label{discr&res} $(R_d^{\deg(\Delta_d)}, \Delta_d^{\deg(R_d)})$
is K-semistable with respect to $SL(d+1,\mathbb{C})$.
\end{claim}
The claim follows from part three of corollary (\ref{examples}) .}
\end{example}
 
K-(semi)stability is formulated in terms of a numerical criterion modeled after Hilbert and Mumford's Geometric invariant theory. We make this explicit by introducing the following
\begin{definition}
\emph{Let $V$ be a rational representation of $G$, and let $\lambda$ be any degeneration in $H$(a maximal algebraic torus of $G$) . The \textbf{\emph{weight}}  $w_{\lambda}(v)$  of $\lambda$ on $v\in V\setminus\{0\}$ is the integer}
\begin{align*}
w_{\lambda}(v):= \mbox{\emph{min}}_{ \{ x\in \mathcal{N}(v)\}}\ l_{\lambda}(x)= \mbox{\emph{min}} \{ <\chi,\lambda>| \chi \in \mbox{\emph{supp}}(v)\}\ .
\end{align*}
 \noindent\emph{ Alternatively, $w_{\lambda}(v)$ is the unique integer such that}
\begin{align*}
\lim_{|t|\rightarrow 0}t^{-w_{\lambda}(v)}\lambda(t)v \  \mbox{ {exists in $V$ and is \textbf{not} zero}}.
\end{align*}
\end{definition}
The precise relationship between weights and K-(semi)stability is brought out in the following
\begin{proposition}\label{weightinequality}
$(v,w)$ is K-semistable if and only if 
\begin{align}
w_{\lambda}(w)\leq w_{\lambda}(v)
\end{align}
for all degenerations $\lambda$ in $G$.
\end{proposition}
Next we equip $V$ and $W$ with \emph{Hermitian norms} which we denote by $||\ ||$ . 
Observe that for $v\in V\setminus\{0\}$ (for example)  we have the following asymptotic expansion
\begin{align} \label{asymp}
\lim_{|t| \rightarrow 0}\log ||\lambda(t)v ||^{2}=w_{\lambda}(v )\log|t|^2+O(1)\ .
\end{align}

\begin{definition}The \textbf{energy} of the pair $(v ,w)$ is the function on $G$ given by
\begin{align}\label{kempfness}
p_{w,v }(\sigma):= \log\frac{||\sigma \cdot w ||^2}{||w ||^2}-\log\frac{||\sigma \cdot v ||^2}{|| v ||^2}\quad \sigma\in G \ .
\end{align}
\end{definition} 
It follows at once from (\ref{asymp}) that the asymptotic behavior of the energy of the pair along any degeneration $\lambda$ is given by
\begin{align}\label{enasymp}
  p_{w,v}(\lambda(t))=\left(w_{\lambda}(w)-w_{\lambda}(v)\right)\log|t|^2+O(1)\ .
\end{align} 
 \begin{definition}\emph {Let $G$ be a reductive algebraic group over $\mathbb{C}$. Consider rational representations $\rho_{V}:G\ra GL(V)$  (respectively $\rho_{W}:G\ra GL(W)$) .  Then $G$ has \textbf{\emph{property P}} if and only if for any pair $(v  ,w)\in V\setminus\{0\}\oplus W\setminus\{0\}$  the following statements are equivalent}
\begin{align*}
\begin{split}
& i)\ \mbox{There exists a degeneration $\lambda$ such that}\ \lim_{|\alpha|\ra 0}  p_{w,v}(\lambda(\alpha)) =-\infty \ . \\
\ \\
& ii)\ \mbox{ There is  a sequence }\   \{\sigma_j\}\subset  G  \ \mbox{such that} \ \lim_{j\ra \infty}  p_{w,v}(\sigma_j)=-\infty \ .
 \end{split}
\end{align*}
 
\end{definition}
\begin{proposition}[Sun's lemma]\label{sun'slemma}
Let $G\cong (\mathbb{C}^*)^N$ be an algebraic torus . Then $G$ has property P.
\end{proposition}
 
 K-semistability, the energy function, and property P are related as follows.
\begin{proposition}\label{kempfnessP} Let $G$ be a reductive algebraic group.
Fix $(v,w)$. Then the following are equivalent 
\begin{align}
\begin{split}
& a)\  p_{w,v}  \ \mbox{is bounded below along all degenerations}\ . \\
\ \\
& b)\  p_{w,v} \ \mbox{is bounded below along all algebraic tori .} \\
\ \\
& c)\ (v,w) \ \mbox{is K-semistable .} \\
 \ \\
&\mbox{If $G$ has property P, then any one of the above implies that }\\
\ \\
& d)\ \mbox{$ p_{w,v}$ is bounded below on $G$.}
\end{split}
\end{align}
\end{proposition}
 \begin{remark}
\emph{In particular, whether or not $ p_{w,v}$ is bounded from below depends only on the pair $(v ,w)$ and not on the norms .}
 \end{remark}
\begin{definition}
\emph{$ p_{w,v}$ is \textbf{\emph{proper}} on $S\subseteq G$ if and only if for all $m>>0$
\begin{align}
 p_{w,v}(\sigma)+ \frac{1}{m}\log\frac{||\sigma \cdot v||^2}{||v||^2}\geq \frac{q(V)}{m}\log||\sigma||_{HS}^2-B \ , \mbox{  $\sigma\in S$}\ .
\end{align}
  $B=B(v,w, ||\ ||, S) $ is a positive constant, and $||\sigma||_{HS}$ denotes the Hilbert Schmidt norm of $\sigma $ with respect to some Hermitian metric on the standard representation.}
\end{definition}  
Applications to K\"ahler Einstein manifolds with discrete automorphism groups require the following proposition.
\begin{proposition}\label{kempfnessproper} The following statements are equivalent.\\
a) The pair $(v,w)$ is K-stable in the strict sense.\\
\ \\
b) $ p_{w,v}$ is proper along all degenerations $\lambda$ in $G$. \\
\ \\
c) $ p_{w,v}$ is proper along all algebraic tori $H\leq G$.
\end{proposition}

%%%%%%%%%%%%%%%%%%%%%%%%%%%%%%%%%%%%%%%%%%%%%%%%%%%%%%%%%%%%%%%%%%%%%%%%%%%
 \subsection{K-stability of complex projective varieties} 
 
 A nontrivial special case of K-stability arises in connection with complex projective varieties .  In order to proceed, let us first recall the Hilbert-Mumford stability theory. The core of this theory consists in associating to a vector bundle $\mathcal{E}$ over a curve $X$ (for example) or a subvariety $X\ra \cpn$ a ``projective geometric gadget" that encodes the object up to projective equivalence. More precisely one associates to these data an \emph{orbit} $G\cdot v$ of some nonzero vector $v$ in a finite dimensional complex rational $G$ module $E$. For example to $\mathcal{E}\ra X$ one associates the \emph{Gieseker point} and to a subvariety $X\ra \cpn$ one associates either the \emph{Hilbert point} or the \emph{Chow form} . Similarly, in order to apply K-stability to a smooth projective variety $X\ra \cpn$ we must associate to our embedded variety $X$ a \emph{pair} $(v_1(X),v_2(X))$ ,  $v_i(X)\in E_i\setminus\{0\}$ where each $E_i$ is a finite dimensional rational $G$-representation. The notation is intended to suggest that $X$ is ``encoded" by the pair $(v_1,v_2)$. As the reader shall see, each $v_i$ is \emph{projectively natural} and by this we mean 
\begin{align}
\sigma\cdot v_i(X)=v_i(\sigma X) \quad \mbox{for all $\sigma\in G$}\ .
\end{align}

\begin{definition}[\textbf{\emph{Cayley-Chow Forms}}]
  \emph{Let $X^n \hookrightarrow \cpn$ be an irreducible, linearly normal subvariety of degree $d$ . The Cayley-Chow form of $X$, denoted by $R_X$,  is the   defining polynomial (unique up to scaling ) of the divisor  }
\begin{align}
\{L\in \gr |\ L\cap X\neq \emptyset \} = \{L\ |\ R_X(L)=0\} \ .
\end{align}
\emph{$R_X$ has degree $d$ in the Pl\"ucker coordinates, moreover the irreducibility of $X$ implies that $R_X$ is also irreducible . }
\end{definition} 
\begin{definition} \emph{ Let $X^n \hookrightarrow \cpn$ be a nonlinear, linearly normal subvariety of degree $d$ . The \textbf{\emph{dual variety}} to $X$, denoted by $X^{\vee}$, is the variety of tangent hyperplanes to $X$}
 \begin{align}
X^{\vee}= \mbox{Zariski closure}\big(\{ f\in {\cpn}^{\vee}|\  \mathbb{T}_p(X)\subset \mbox{ker}(f) \quad \mbox{for some $p\in X\setminus X_{sing}$} \}\big) \ .
\end{align}
\end{definition}
$ \mathbb{T}_p(X)$ denotes the \emph{embedded} tangent space to $X$ at the point $p$. $ \mathbb{T}_p(X)$ is an $n$ dimensional linear subspace of $\cpn$.
 
 \begin{definition}
 \emph{The \textbf{\emph{dual defect}} of $ X\hookrightarrow \cpn$ is the nonnegative integer 
\begin{align} 
\delta(X):=N-\mbox{dim}(X^{\vee})-1
\end{align}}
\end{definition}
Most varieties have dual defect equal to 0 . There is a well known upper bound on the defect which we require
for the definition of K-stability .\\
\ \\
\noindent\textbf{Theorem .} (F.L. Zak \cite{zak}) \emph{Let $X^n\ra \cpn$ ($n\geq 2$)  be a linearly normal irreducible variety which is not a linear space.   Then}
\begin{align}\label{zak}
\delta\leq n-2 \ .
\end{align}

When $X^{\vee}$ is indeed a hypersurface (i.e. $\delta=0$) following Gelfand, Kapranov and Zelevinsky the defining polynomial, unique modulo scaling, is denoted by $\Delta_X$, which we shall call the \emph{\textbf{${X}$-discriminant}.}
\begin{align}
X^{\vee}=\{ f\in {\cpn}^{\vee}|\ \Delta_X(f)=0\} \ .
\end{align}
 
Just as in the case of resultants and discriminants of polynomials in one variable, we may view the general $X$-discriminant and Cayley-Chow form as homogeneous polynomials on spaces of matrices:
\begin{align}
\begin{split}
&\Delta_{X}\in \mathbb{C}[M_{1\times(N+1)}  ] \\
\ \\
&R_X\in \mathbb{C}[M_{(n+1)\times(N+1)}] \ . \\
\end{split}
\end{align}
The action of $\sigma \in GL(N+1,\mathbb{C})$ on these two polynomials is given by
\begin{align}\label{action}
\begin{split}
&\sigma\cdot \Delta_{X}((a_{i}))=\Delta_{X}( (a_{i})\cdot\sigma) \\
\ \\
&\sigma \cdot R_X((c_{ij}))=R_X((c_{ij})\cdot \sigma) \ .
\end{split}
\end{align}
Next, fix $k\in \mathbb{N}_+$ and $(l_1,l_2,\dots,l_k)$ with $l_i\in \mathbb{N}_+$ . We set $\mathbb{P}^{(l_{\bull})}:=\mathbb{P}^{l_1}\times\dots\times\mathbb{P}^{l_k}$ .
Consider the Segre embedding
\begin{align*}
X\times\mathbb{P}^{(l_{\bull})}\ra \mathbb{P}(\mathbb{C}^{N+1}\otimes\mathbb{C}^{(l_1+1)}\otimes\dots\otimes\mathbb{C}^{(l_k+1)})
\end{align*}
\begin{definition}
Assume the dual defect of $X\times\mathbb{P}^{(l_{\bull})}$ vanishes. The \textbf{$X$-hyperdiscriminant} of format $(l_{\bull})$ is the irreducible defining polynomial $\Delta_{(l_{\bull})}$ of $(X\times\mathbb{P}^{(l_{\bull})})^{\vee}$.
\end{definition}
The hyperdiscriminant $\Delta_{(l_{\bull})}$ is an irreducible polynomial in the entries of a``hypermatrix"
\begin{align*}
\Delta_{(l_{\bull})}\in \mathbb{C}[M_{(l_1+1)\times\dots\times(l_k+1)\times(N+1)}] \ .
\end{align*}

The circumstances which insure that the  {dual defect} of the Segre image  of $X\times\mathbb{P}^{(l_{\bull})} $ is equal to zero has been completely worked out by Weyman and Zelevinsky in \cite{caytrick}. In these cases we say that the hyperdiscriminant is well-formed. When $N=n$ and therefore $X=\mathbb{P}^n$
the hyperdiscriminant is the \emph{hyperdeterminant} of Cayley, Gelfand, Kapranov, and Zelevinsky see \cite{hyperdet}.  {What is relevant for my applications to the Mabuchi energy are the hyperdiscriminants of format $(n-1)$}.\\
\ \\
\noindent\textbf{Theorem} (Weyman, Zelevinsky \cite{caytrick})\textbf{.}
\emph{Let $X^n$ be an $n$ dimensional, linearly normal subvariety of $\cpn$ where $N>n$, then
 the X-hyperdiscriminant of format $(l_{\bull})$ exists if and only if the following two inequalities hold}
\begin{align}\label{trick}
\begin{split}
& a)\ l_i\leq n+\sum_{i\neq j}l_j \quad 1 \leq i\leq n \\
\ \\
&  b)\ \delta\leq  \sum_{1\leq i\leq k}l_i  \ . \\
\ \\
&\mbox{\emph{In particular}}\
X\times\mathbb{P}^n,\ X\times\mathbb{P}^{n-1},\dots,X\times\mathbb{P}^{\delta(X)}
\mbox{ \emph{are all dually nondegenerate}}\\
&\mbox{\emph{in their Segre embeddings. Moreover,}}\\
 \ \\
&i) \ \Delta_{X\times\mathbb{P}^n}=R_X \quad (\mbox{\emph{the ``Cayley trick" }}) \\
\ \\
&ii)\ \Delta_{X\times\mathbb{P}^{\delta(X)}}=R_{X^{\vee}}\quad (\mbox{\emph{the ``dual Cayley trick" }})\ .
\end{split}
\end{align}
 
When $X$ is a \emph{smooth} subvariety we may make use of a result due to Beltrametti, Fania, and Sommese which  exhibits the degree and codimension of the dual in terms of the top Chern class of the \emph{jet bundle} $J_1(\mathcal{O}_X(1))$ (see \ref{1jet}). This result is used extensively in the main argument of the paper, we shall use it to find the degree of the hyperdiscriminant. \\
\ \\
\textbf{Theorem} (Beltrametti, Fania, and Sommese \cite{bfs})\textbf{.} 
\emph{Assume $X$ is smooth. Then $X^{\vee}$ is a hypersurface if and only if $c_n(J_1(\mathcal{O}_X(1)))\neq 0$. Moreover,}\\
\ \\
$ i)\ \mbox{deg}(\Delta_X)=\int_{X}c_n(J_1(\mathcal{O}_X(1))) $. \\
\ \\
\emph{ More generally, when $\delta(X)>0$, we have the following}\\
\ \\
$ii)\ \mbox{{deg}}(X^{\vee})=\int_{X}c_{n-\delta(X)}(J_1(\mathcal{O}_X(1)))\om^{\delta(X)}$ .\\
\ \\
$ iii)\ \delta(X)=\min\{k \ | \ c_{n-k}(J_1(\mathcal{O}_X(1)))\neq 0\}$ .\\
\ \\
 
\begin{definition} \emph{Let $X\hookrightarrow \cpn$ be a linearly normal $n$ dimensional variety with degree $d\geq 2$. Fix a maximal algebraic torus $H\leq G$.
The weight polytope of the $X$-resultant $\mathcal{N}(R_X)$ is called the \emph{\textbf{Chow polytope}} of $X$  and the weight polytope of the $X$-hyperdiscriminant $\mathcal{N}(\Delta_{X\times\mathbb{P}^{n-1}})$ is the \emph{\textbf{hyperdiscriminant polytope}}.}
\end{definition}

\begin{remark}
\emph{Once more, the reader should bear in mind that there are smooth varieties $X$ whose dimensions exceed 2 such that $\delta(X)>0$ ( for example $\mathbb{P}^2\times \mathbb{P}^{1}$ in its Segre embedding ) . 
Zaks' bound  $\delta(X)\leq n-2$ (see (\ref{zak}) ) implies that the hyperdiscriminant $\Delta_{\xhyp}$ is well formed \footnote{The reader should realize that this amounts to the fact that $\Delta_{\xhyp}$ is a \emph{non-constant} polynomial.} for \emph{any} $X$ (irreducible, linearly normal, $\mbox{deg}(X)\geq 2$) .}
\end{remark}

 With these preparations we introduce the following new stability concept for complex projective varieties, which we call \emph{K-(semi)stability} in order to avoid proliferation of terminology.  
Our new idea extends the concept of K-semistability proposed by Tian in \cite{psc}. Our definition seems to be quite different from that proposed by Donaldson in \cite{skdtoric} and developed by his many followers.  The reader should keep in mind that the crucial difference between our formulation and the conventional one is that from our new viewpoint \emph{the limit cycle plays no role.}
\begin{definition}\label{fundamental}  \emph{Let $X\ra \cpn$ be a nonlinear, linearly normal, complex projective variety (not necessarily smooth). Then $X$ is \emph{\textbf{K-semistable}} provided the pair of polynomials
\begin{align*}
\left(R_X^{\deg(\Delta_{X\times\mathbb{P}^{n-1}})},\ \Delta_{X\times\mathbb{P}^{n-1}}^{\deg(R_X)}\right)  \qquad (*)
\end{align*}
is K-semistable in the sense of definition \ref{kstability} with respect to the natural action of $G=SL(N+1,\mathbb{C})$ on the irreducible modules
\begin{align}
\begin{split}
&\Delta_{X\times\mathbb{P}^{n-1}}^{\deg(R_X)}\in \mathbb{C}_{\deg(\Delta)\deg(R)}[M_{n\times(N+1)}  ]^{SL(n,\mathbb{C})} \\
\ \\
& R_X^{\deg(\Delta_{X\times\mathbb{P}^{n-1}})}\in \mathbb{C}_{\deg(\Delta)\deg(R)}[M_{(n+1)\times(N+1)}]^{SL(n+1 ,\mathbb{C})}\ .
\end{split}
\end{align}
  That is, for all maximal algebraic tori $H\leq G$ the scaled hyperdiscriminant polytope dominates the scaled Chow polytope}
 \begin{align}\label{include1}
\deg(\Delta_{X\times\mathbb{P}^{n-1}} )\mathcal{N}(R_X)\subseteq  {\deg(R_X)}\mathcal{N}(\Delta_{X\times\mathbb{P}^{n-1}}) \ .
\end{align}
\emph{$X$ is \emph{\textbf{K-stable}} (in the strict sense) if and only if the pair $(*)$  is K-stable where the degree $q$ is given by $q=\deg(\Delta_{X\times\mathbb{P}^{n-1}} )\deg(R_X) $ .}
\emph{A polarized algebraic variety $(X,L)$ is \emph{\textbf{asymptotically K-(semi)stable}} provided $X\overset{L^r}{\ra}\mathbb{P}^{N_r}$ is K-(semi)stable for all $r>>0$ where $m_0$ is independent of $r$ .}
\end{definition}
In the definition we have abused notation by setting $\deg(\Delta):= \deg(\Delta_{X\times\mathbb{P}^{n-1}})$.
\begin{remark}
\emph{$H$ must be allowed to vary in our definition. This is due to our requirement that the K stability of $X$ imply (and be implied by) the K stability of any subvariety of $\cpn$ \emph{projectively equivalent} to $X$.}
\end{remark}
%%%%%%%%%%%%%% %%%%%%%%%%%%%%
 Theorem B, together with the considerations of the previous section on the energy asymptotics of pairs (see (\ref{enasymp})) completely justify the following definition.
\begin{definition} \emph{ Let $X$ be a smooth, linearly normal subvariety of  $\cpn$.
Fix any maximal torus $H$ of $G$ and let $\lambda$ be any degeneration in $H$. Then the \textbf{\emph{generalized Futaki invariant}} $F_P(\lambda)$ of $\lambda$ is the integer given by}
\begin{align}
F_P(\lambda):={\deg(R_X)}\mbox{\emph{min}}_{ \{ x\in \mathcal{N}(\Delta_{X\times\mathbb{P}^{n-1}})\}}\ l_{\lambda}(x)- 
 {\deg(\Delta_{X\times\mathbb{P}^{n-1}} )} \mbox{\emph{min}}_{ \{ x\in \mathcal{N}(R_X)\}}\ l_{\lambda}(x) \ .
 \end{align}
\end{definition}
  The following is a special case of proposition (\ref{weightinequality}) .
\begin{proposition}
$X\hookrightarrow \cpn$ is K-semistable if and only if the generalized Futaki invariant is less than or equal to zero for all degenerations $\lambda$ in $G$.
\end{proposition}

To close this section, we need to discuss the relationship between our encoding process and \emph{limit cycle formation}.  
As we have mentioned, in Mumford's G.I.T. $v(X)$ may be given in terms of Hilbert points or Chow forms. 
In both cases it is known that the encoding is natural
\begin{align}
\mbox{Hilb}_m(\sigma\cdot X)=\sigma\cdot \mbox{Hilb}_m(X) \ , \ R_{\sigma\cdot X}=\sigma\cdot R_X \ .
\end{align}
Let $\lambda$ be an algebraic one parameter subgroup of $G=\slnc$, for a given $X\subset \cpn$ we let ${\lambda(0)}X$ denote the flat limit cycle of $X$ under $\lambda$. This is considered to be a point in the Hilbert scheme. In Mumford's theory, a crucial property of Hilbert points and Cayley -Chow forms is the following compatibility with cycle formation
\begin{align}
\mbox{Hilb}_m({\lambda(0)}X)=\lambda(0)\cdot \mbox{Hilb}_m(X) \ , \ R_{{\lambda(0)}X}=\lambda(0)\cdot R_X \ .
\end{align}
In K-stability this compatibility \emph{fails}. Simply put, ${\lambda(0)}X$ in general has no meaningful tangent plane and therefore \emph{the hyperdiscriminant is undefined} .   
 %%%%%%%%%%%%%%%%%%%%%%%%BOTTCHERN%%%%%%%%%%%%%%%%%%%%%%%%%%
\section{Bott-Chern Forms and Donaldson Functionals}
Let $\phi$ be a $GL_n(\mathbb{C})$ invariant polynomial on $M_{n\times n}(\mathbb{C})$ homogeneous of degree $d$. The \emph{complete polarization} of $\phi$ is defined as follows. Let $\tau_1,\tau_2,\dots, \tau_d$ be arbitrary real parameters. Then
\begin{align*}
\phi(\tau_1A_1+\tau_2A_2+\dots +\tau_dA_d)=\sum_{|\alpha|=d}\phi_{\alpha}(A_1,A_2,\dots,A_d)\tau^{\alpha} \  \quad \tau^{\alpha}:=\tau_1^{\alpha_1}\tau_2^{\alpha_2}\dots\tau_d^{\alpha_d} \ .
\end{align*}
We let $\phi_{(1)}(A_1,A_2,\dots,A_d)$ denote the coefficient of $\tau_1\tau_2\dots\tau_d$ . We define
\begin{align*}
\phi_{(1)}(A;B):=\phi_{(1)}(A,\overbrace{B,B,\dots,B}^{d-1}) \ .
\end{align*}
 Let $M$ be an $n$ dimensional complex manifold, $E$ is a holomorphic vector bundle of rank $k$ over $M$. $H_0$ and $H_1$ are two Hermitian metrics on $E$.
Let $H_t$ be a smooth path joining $H_0$ and $H_1$ in $\mathcal{M}_{E}$ ( the space of Hermitian metrics on $E$) . Define $U_t:= (\frac{\dl}{\dl t}H_t)\cdot H^{-1}_t $ . $F_t:= \dlb\{(\dl H_t)H^{-1}_t\}$ is the curvature of $H_t$ (a purely imaginary (1,1) form) .  Now suppose that $\phi$ is a homogeneous invariant polynomial on $M_{k\times k}(\mathbb{C})$ of degree $d$. Then 
\begin{align}
\phi_{(1)}(U_t;\ F_t)
\end{align}
is a form of type $(d-1, d-1)$ . Observe that the following identity holds, this is used below.
\begin{align}
\deg(\phi)\phi_{(1)}(A\ ;\ B)=\frac{\dl}{\dl b}\phi(B+bA)|_{b=0} \ .
\end{align}

The Bott Chern form is given as follows 
\begin{align}
BC(E,\phi ; H_0,H_1):=-\frac{\deg(\phi)}{(n+1)!}\int_0^1\phi_{(1)}(U_t;\ \frac{\sqrt{-1}}{2\pi}F_t)\ dt \ .
\end{align}
\begin{proposition}[R. Bott , S.S. Chern \cite{bottchern}]
\begin{align}
\frac{\sqrt{-1}}{2\pi} \dl\dlb BC(E,\phi ; H_0,H_1)=\phi(\frac{\sqrt{-1}}{2\pi}F_1)-\phi(\frac{\sqrt{-1}}{2\pi}F_0) \ .
\end{align} 
\end{proposition}
When ${\deg}(\phi)$ has degree $n+1$ $BC(E,\phi ; H_0,H_1)$ is a \emph{top dimensional} form on $M$, and the following integral is well defined
\begin{align}
D_{E}(\phi; H_0,H_1):= \int_M BC(E,\phi ; H_0,H_1) \ .
\end{align}
Of particular importance is when $\phi (A)=Ch_{n+1}(A)=\frac{1}{(n+1)!}\mbox{Tr}(A^{n+1})$ . Observe that in this case we have
\begin{align}
 \phi_{(1)}(A;B)=\mbox{Tr}(AB^n) \ .
\end{align}
Let $H:Y\rightarrow \mathcal{M}_{{E}}$ (the space of $C^{\infty}$ Hermitian metrics on $E$)  be a smooth map, where $Y$ is a complex manifold of dimension $m$. \emph{Fix} a Hermitian metric $H_0$ on $E$.
Then we are interested in the smooth function on $Y$
\begin{align}
Y\ni y \rightarrow D_{E}(\phi; H_0, H(y)) \ .
\end{align}

Let $p_2$ denote the projection from $Y\times M$ onto $M$. Then $H(y)$ is a smooth Hermitian metric on $p^{*}_2(E)$  
 whose curvature is given by
\begin{align}
F_{Y\times M}(H(y))=\dlb_{Y\times M}\{(\dl_{Y\times M}H(y))H(y)^{-1}\} \ .
\end{align}
For the proof of the following proposition, see \cite{bottchrnfrms} prop. 1.4 on pg. 213 .
\begin{proposition} Let $\phi$ be homogeneous of degree $n+1$ and $H_0$ a fixed metric on ${E}$. Then for all smooth compactly supported forms $\eta$ of type $(m-1,m-1)$ we have the identity
\begin{align}\label{workhorse}
\frac{\sqrt{-1}}{2\pi}\int_YD_{E}(\phi; H_0, H(y))\dl_Y\dlb_Y\eta=\int_{Y\times M}\phi(F_{Y\times X}(\frac{\sqrt{-1}}{2\pi}H(y)))\wedge p_1^*(\eta)\ .
\end{align}
\end{proposition}

Now we extend the above to holomorphic Hermitian \emph{complexes} $(E^{\bull}, H_0^{\bull} ; \dl^{\bull})$ .
Let $H^{\bull}:Y\rightarrow \mathcal{M}_{E^{\bull}}$ be a smooth map. Concretely, $H^{\bull}(y)$ is a $C^{\infty}$ metric on $E^{\bull}$. We define the smooth function on $Y$
\begin{align}
D_{E^{\bull}}(\mbox{Ch}_{n+1}; H_0^{\bull}, H^{\bull}(y)):=
\sum_{j=0}^{l}(-1)^jD_{E^{j}}(\mbox{Ch}_{n+1}; H_0^{j}, H^{j}(y)) \ .
\end{align}

\begin{corollary}\label{ddbardfnctl}
For all smooth compactly supported $(m-1,m-1)$ forms $\eta$ on $Y$ we have
\begin{align}
\begin{split}
&\int_Y\frac{\sqrt{-1}}{2\pi}D_{E^{\bull}}(\mbox{Ch}_{n+1}; H_0^{\bull}, H^{\bull}(y))\wedge \dl_Y\dlb_Y \eta=\\ 
\ \\
&\sum_{j=0}^{l}(-1)^j\int_{Y\times M}\mbox{Ch}_{n+1}(\frac{\sqrt{-1}}{2\pi}F_{Y\times X}^{E^j}(H^j(y)))\wedge p_1^*(\eta) \ .
\end{split}
\end{align}
\end{corollary}

%%%%%%%%%%%%%%%%%%%%%%%%%%%%%%%%%%%%%%%%%%%%%%%%%%%%%%%%%%%%%%%%%%%%%%%%%%%%%%%%%% 
 
  \section{The Main Lemma}
 Let $X\rightarrow \cpn$ be a smooth, linearly normal, $n$-dimensional, and dually non-degenerate complex projective variety.  Let  $\mathbb{G}(n,\cpn)$ denote the Grassmannian of $n$ dimensional linear subspaces of $\cpn$ and let $\rho:X\ra \mathbb{G}(n,\cpn)$ denote the \emph{Gauss map} of $X$
 \begin{align}
 \rho(p):=\mathbb{T}_p(X)\in \mathbb{G}(n,\cpn) 
 \end{align}
 where $\mathbb{T}_p(X)$ denotes the \emph{embedded tangent space} to $X$ at $p$. Let $\mathcal{U}$ denote the rank $n+1$ universal (tautological) bundle over $\mathbb{G}(n,\cpn)$.
   Of basic importance throughout the paper is the pull back of this bundle under the Gauss map 
 $\rho^*(\mathcal{U})$ which we shall denote by $J_1(\mathcal{O}_X(1))^{\vee}$.  This is (dual to) the bundle of \emph{one jets} of $\mathcal{O}(1)|_X$ .  We apply the construction of the previous section to the terms of a \emph{Cayley-Koszul complex} on $X$, defined as follows
 
 \begin{align}
&K^i :=\bigwedge^{i}J_1(\mathcal{O}_X(1))^{\vee}  \ .
\end{align}
 
 It is clear that $J_1(\mathcal{O}_X(1))^{\vee}$ is a subbundle of the trivial bundle $X\times \mathbb{C}^{N+1}$ and so inherits the standard euclidean (Hermitian ) metric from $\mathbb{C}^{N+1}$. In this way, as in the previous section, we may define maps $H^{\bull}:G\ra \mathcal{M}_{K^{\bull}}$ , where $G=\slnc$ plays the role of $Y$.
 
\begin{align}
&H^i(\sigma)=\bigwedge ^{i}(h_{\cn} \circ \sigma)|_{J_1(\mathcal{O}_X(1))^{\vee}}  \ .
\end{align}

$h_{\cn}$ denotes the standard Hermitian form on $\cn$
\[h_{\cn}(V,W):= v_0\bar{w_0}+v_1\bar{w_1}+\dots +v_N\bar{w_N} \ . \]

 \noindent\textbf{Main Lemma .}
\emph{ Let $X\hookrightarrow \cpn$ be a smooth, linearly normal $n$ dimensional subvariety. Let $X^{\vee}$ be the dual of $X$. Assume that $X^{\vee}$ is a hypersurface with defining polynomial  $\Delta_X$.  Then  there is a continuous norm $||\ ||$ on the vector space of degree $d^{\vee}:=\deg(X^{\vee})$ polynomials on $(\mathbb{C}^{N+1})^{\vee}$  such that 
 \begin{align}\label{mainlemma}
 D_{\ K^{\bull} }(Ch_{n+1} ; H^{\bull}(\sigma),\ H^{\bull}(e) \ ) =\log\frac{||\sigma \cdot \Delta_X||^2}{||\Delta_X||^2} \ .
 \end{align}
  $\sigma \in G$ ($e$ denotes the identity in $G$) . } \\
\begin{remark}
\emph{The Main lemma exhibits the ``height" of the defining equation of $Z=X^{\vee}$ ( a global, purely algebro geometric object) as an integral over $X$ of a local curvature quantity  derived from the metric $\oms$ }.\\
\end{remark}
 
 Before we proceed to the proof of the main lemma, let us explain what is meant by a continuous metric (or norm) on $\mathcal{O}_{B}(-1)$,  where $B:=\mathbb{P}(H^0({\cpn}^{\vee}, \mathcal{O}({d}^{\vee})))$ and ${d}^{\vee} $ denotes the degree of $X^{\vee}$ .  \emph{Up to scaling} we have that $\Delta_X \in H^{0}({\mathbb{P}^{N}}^{\vee},\mathcal{O}( {d}^{\vee}))$.  

In general we write linear form $f$ on $\cpn$  (i.e. a \emph{point} in the \emph{dual} $\cpn$) as $f= a_0z_0+a_1z_1+\dots +a_Nz_N$. Therefore we take $[a_0:a_1:\dots :a_N]$ as the homogeneous coordinates of $f$ on ${\cpn}^{\vee}$.
Therefore we may write
\begin{align*}
\Delta_X(f)=\sum_{|\alpha|=d^{\vee}}c_{\alpha_0,\dots \alpha_N}{a_0}^{\alpha_0}{a_1}^{\alpha_1}\dots {a_N}^{\alpha_N} \ .
\end{align*}
The \emph{finite dimensional complex vector space} $H^{0}({\mathbb{P}^{N}}^{\vee},\mathcal{O}(d^{\vee}))$ comes equipped  with its standard Hermitian inner product $< ,\ >$ in which the monomials ${a_0}^{\alpha_0}{a_1}^{\alpha_1}\dots {a_N}^{\alpha_N}$ form an orthogonal basis. Under a suitable normalization we have that
\begin{align*}
||\Delta_X||_{FS}^2:=<\Delta_X,\Delta_X>=\sum_{|\alpha|=d^{\vee}}\frac{|c_{\alpha_0,\dots \alpha_N}|^2}{\alpha_0!\alpha_1!\dots \alpha_N!}\ .
\end{align*}
Finally, to say that the metric $||\ ||$ on $\mathcal{O}_{B}(-1)$ is \emph{continuous} means that there is a continuous function $\theta$ on $B$ such that
\begin{align}
\exp({\theta})||\ ||_{FS}= ||\ || \ .
\end{align}
Since $B$ is compact, the conformal factor $\exp({\theta})$ is \emph{bounded}. This is the key point.

We first construct the norm appearing on the right hand side of  (\ref{mainlemma}). Recall that the \emph{universal hypersurface associated to $B$} is given by
\begin{align}
\Sigma:= \{([F],\ [a_0:a_1:\dots :a_N])\in B\times  {\cpn}^{\vee} | \ F(a_0,a_1,\dots, a_N)=0 \} \ .
\end{align}
Then $\Sigma$ is the base locus of the natural section
\begin{align}
\varphi \in H^0(B\times{\cpn}^{\vee},p_1^*\mathcal{O}_B(1)\otimes p_2^*\mathcal{O}_{\cpn}({d}^{\vee})) \quad ,\ \Sigma= \{\varphi=0\} \  .
\end{align}
Let $\om$ denote the K\"ahler form on the dual $\cpn$. We consider the $(1,1)$ current $u$ on $B$ defined by the fiber integral ${p_1}_*p_2^*(\om^{N})$.
\begin{diagram}
  \Sigma &\rTo^{p_2}&{\cpn}^{\vee}\\
 \dTo^{p_1}\\
   B
\end{diagram}
 That is, for all $C^{\infty}$ $(b-1, b-1)$ forms $\psi$ on $B$ we have
\begin{align}\label{current}
\int_{B}u\wedge \psi=\int_{\Sigma}p_2^*(\om^N)\wedge p_1^*(\psi) \ .
\end{align}
For the following, see \cite{psc} Lemma 8.7 pg. 32 .
\begin{proposition} 
The cohomology class of the current $u$ coincides with the class of $\om_{B}$ (the Fubini-Study form). Moreover, there is a continuous function $\theta$ on $B$ such that, in the sense of currents we have
\begin{align}
u=\om_{B}+\frac{\sqrt{-1}}{2\pi}\dl\dlb \theta \ .
\end{align}
\end{proposition}

Let $\mathcal{I}(\sigma):=D_{\ K^{\bull} }(Ch_{n+1} ; H^{\bull}(\sigma),\ H^{\bull}(e) \ )$. The main point is to establish the following proposition.
\begin{proposition}\label{pluri} Let $|| \ \cdot ||:= \exp(\theta)||\ \cdot ||_{FS}$. Then the difference
\begin{align}
\mathcal{I}(\sigma)-  \log  \frac{{||\sigma\cdot\Delta_{X}||}^{2}}{{||\Delta_{X}||}^{2}} 
\end{align}
is a pluriharmonic function on $G$ .
\end{proposition}

We require the simple flag variety
\begin{align}
I_{\Delta}:=\{(L,f)\in \mathbb{G}(n,\cpn)\times {\cpn}^{\vee}\ | \ L\subset \mbox{ker}(f)\}\ .
\end{align}

 \noindent \emph{Proof of proposition \ref{pluri}}.  
 \begin{lemma} Let $p_i$ denote the projection onto the $ith$ factor of the flag variety $I_{\Delta}$.
\begin{diagram}
  I_{\Delta} &\rTo^{p_2}&{\cpn}^{\vee}\\
 \dTo^{p_1}\\
  \mathbb{G}(n,\ N)
\end{diagram}
Let $\om_{{\cpn}^{\vee}}$ the Fubini Study K\"ahler form on ${\cpn}^{\vee}$. Then we have the following identity of forms on $\mathbb{G}(n,\ N)$
\begin{align}\label{pdual}
{p_1}_*(p_2^*(\om^N_{{\cpn}^{\vee}}))=\sum_{i=0}^{n+1}(-1)^i\mbox{\emph{Ch}}(\bigwedge^i{\mathcal{U}}, \ h_{FS})^{\{n+1,n+1\}}\ .
\end{align}
\end{lemma}

To see this, observe that the left hand side of (\ref{pdual}) is of type $(n+1,n+1)$ and invariant under the action of the unitary group. The latter implies that it must be a polynomial in the forms $c_1(\mathcal{U}^{\vee}),c_2(\mathcal{U}^{\vee}),\dots, c_{n+1}(\mathcal{U}^{\vee})$. Let $\Omega$ be any invariant form on  $\mathbb{G}(n,\ N)$ of type complimentary to ${p_1}_*p_2^*\om^N_{ {\cpn}^{\vee}}$. Then
\begin{align}
\begin{split}
\int_{\mathbb{G}(n,\ N)}{p_1}_*(p_2^*(\om^N_{{\cpn}^{\vee}}))\wedge \Omega&=\int_{I_{\Delta}}p_2^*(\om^N_{{\cpn}^{\vee}})\wedge {p_1}^*(\Omega)\\
\ \\
&=\int_{\mathbb{G}(n,\ N)\times {\cpn}^{\vee}}[I_{\Delta}]\wedge p_2^*(\om^N_{{\cpn}^{\vee}})\wedge p_1^*(\Omega) \ .
\end{split}
\end{align}
Observe that $I_{\Delta}=\{s=0\}$ where $s$ is the section of $p_1^*\mathcal{U}^{\vee}\otimes p_2^*\mathcal{O}_{{\cpn}^{\vee}}(+1)$ given by evaluation. Therefore
\begin{align*}
 [I_{\Delta}]&=c_{n+1}(p_1^*\mathcal{U}^{\vee}\otimes p_2^*\mathcal{O}_{{\cpn}^{\vee}}(+1))\\
\ \\
&= \sum_{i=0}^{n+1}c_{1}(p_2^*\mathcal{O}_{{\cpn}^{\vee}}(+1))^{n+1-i}\wedge c_{i}(p_1^*\mathcal{U}^{\vee}) \\
\ \\
&=c_{n+1}(p_1^*\mathcal{U}^{\vee})+\sum_{i=0}^{n}c_{1}(p_2^*\mathcal{O}_{{\cpn}^{\vee}}(+1))^{n+1-i}\wedge c_{i}(p_1^*\mathcal{U}^{\vee}) \ .
\end{align*}
Therefore, for \emph{all} invariant forms $\Omega$ (of complimentary type) we have
\begin{align*}
\int_{\mathbb{G}(n,\ N)}{p_1}_*(p_2^*(\om^N_{{\cpn}^{\vee}}))\wedge \Omega=\int_{\mathbb{G}(n,\ N)} c_{n+1}(\mathcal{U}^{\vee})\wedge \Omega \ .
\end{align*}
Therefore,
\begin{align*}
{p_1}_*(p_2^*(\om^N_{{\cpn}^{\vee}}))=c_{n+1}(\mathcal{U}^{\vee},\ h_{FS}) \ .
\end{align*}
Then the lemma follows immediately from the well known \emph{Borel-Serre identity}  (see \cite{fulton} pg. 57 example 3.2.5 ) .
\begin{align}
\sum_{j=0}^{k}(-1)^j\mbox{Ch}(\bigwedge^jE^{\vee})=c_{k}(E)\mbox{Td}(E)^{-1} \quad (k=rnk(E)) \ .
\end{align}
 
 Let $GX$ be given by
 \begin{align}
 GX:=\{(\sigma, y)\in G\times\cpn \ |\ y\in \sigma X\ \} \ .
 \end{align}
 
 There is a natural map $\rho_G:GX\ra \mathbb{G}(n,\cpn)$ given by
 \begin{align}
 \rho_G(\sigma,y)=\mathbb{T}_y(\sigma X) \ .
 \end{align}

We let $K^i_G$ denote the complex on $GX$ whose fiber at $(\sigma,y)$ is $K^i_y(\sigma X)$. That is, $K^i_G$ is obtained by pulling back $\bigwedge^{i}\mathcal{U}$ via $\rho_G$.  Then we have the identity of forms on $GX$  .
\begin{align}
\sum_{i=0}^{n+1}(-1)^i\rho_{G}^{*}\mbox{Ch}(\bigwedge^i \mathcal{U},\ h_{FS})=\sum_{i=0}^{n+1}(-1)^i\mbox{Ch}(K^i_{G}, \ \rho_G^{*}(h_{FS})) \ .
\end{align}
 Now we consider the diagram
 \begin{diagram}
 {\rho_{G}}^*(I_{\Delta})  &\rTo^{\pi_2}&I_{\Delta} &\rTo^{p_2}&{\cpn}^{\vee}\\
\dTo ^{\pi_{1}}&&\dTo^{p_1}\\
 GX&\rTo ^{\rho_{G}}& \mathbb{G}(n,\ N)&\\
\dTo^{\pi}\\
G
\end{diagram}
 Let $\eta$ be a smooth compactly supported form on $G$ of type
 $(N^2+2N,N^2+2N)$. An application of corollary (\ref{ddbardfnctl}) gives
 \begin{align*}
 \int_G\frac{\sqrt{-1}}{2\pi}\dl\dlb \mathcal{I}\wedge \eta&=\int_{GX}\sum_{i=0}^{n+1}(-1)^{i}\mbox{Ch}_{n+1}( K^i_{G} , \ \rho_G^{*}(h_{FS}))\wedge \pi^{*}(\eta)\\
\ \\
&=\int_{GX}\rho^*_{G }({p_1}_*(p_2^*(\om^N_{{\cpn}^{\vee}})))\wedge \pi^*(\eta) \\
\ \\
&= \int_{\rho^*_{G }(I_{\Delta})}{\pi_2}^*(p_2^*(\om^N_{{\cpn}^{\vee}}))\wedge \pi_1^*(\pi^*(\eta)) \ .
 \end{align*}
Below $T$ denotes the evaluation map $T(\sigma):=[\sigma\cdot \Delta_{X}]$ and $\Sigma$ denotes the universal hypersurface for the family $B:=\mathbb{P}(H^0({\cpn}^{\vee}, \mathcal{O}(d^{\vee})))$.
\begin{diagram}
  T^*(\Sigma)  &\rTo^{\pi_2}&\Sigma &\rTo^{p_2}&{\cpn}^{\vee}\\
\dTo ^{\pi_{1}}&&\dTo^{p_1}\\
 G&\rTo ^{T}& B&\\
 \end{diagram}
Let $u$ denote the positive current defined in (\ref{current}). Using the notation and commutativity in the diagram above gives that
\begin{align}
\begin{split}
\int_GT^*(u)\wedge \eta &=\int_{T^*(\Sigma)}\pi_2^*(p_2^*(\om^N_{{\cpn}^{\vee}}))\wedge \pi_1^*(\eta)\\
 \ \\
&=\int_{\rho^*_{G }(I_{\Delta})}\pi_2^*(p_2^*(\om^N_{{\cpn}^{\vee}}))\wedge\pi^*_1\pi^*(\eta) \ .
\end{split}
\end{align}
We have used that ${T^*(\Sigma)}\cong \rho^*_{G}(I_{\Delta}) $ (birational equivalence). We remark that this holds only because of our assumption that $X$ is dually non-degenerate.  By definition we have that
\begin{align}
T^*(u)=\frac{\sqrt{-1}}{2\pi}\dl\dlb\log\left(e^{\theta\circ T}\frac{||\sigma\cdot\Delta_X||_{FS}^2}{||\Delta_X||_{FS}^2}\right) \ .
\end{align}
Therefore,
\begin{align}
\begin{split}
&\int_G\dl\dlb\left(\mathcal{I}(\sigma)-\log\left(e^{\theta\circ T}\frac{||\sigma\cdot\Delta_X||_{FS}^2}{||\Delta_X||_{FS}^2}\right)\right)\wedge \eta=0\ .\\
\end{split}
\end{align}
For all compactly supported forms $\eta$. Hence the difference is pluriharmonic.
This establishes Proposition \ref{pluri} .

 Since $G$ is simply connected there is an entire function $F$ on $G$ such that
\begin{align}
\mathcal{I}(\sigma)-\log\left(e^{\theta\circ T}\frac{||\sigma\cdot\Delta_X||_{FS}^2}{||\Delta_X||_{FS}^2}\right)
=\log(|F(\sigma)|^2) \ .
\end{align}
 The argument from \cite{psc} (see Lemma 8.8 pg. 34) shows that $F\equiv 1$. This completes the proof of the main lemma. $\Box$

In higher dimensions ($n\geq 2$) two problems with (\ref{mainlemma}) arise which reveal that the $X$-discriminant (the right hand side of (\ref{mainlemma})) must be \emph{modified} in order to make contact with the Mabuchi energy. The first problem we encounter is that when $n\geq 3$ \emph{the projective dual to $X$  may fail to have codimension one}, e.g. take $\mathbb{P}^1\times\mathbb{P}^2$ in its Segre embedding. In a situation like this $\Delta_X$ is taken to be a conveniently chosen constant. The second problem in higher dimensions ($n\geq 2$) is that the global Donaldson (energy) functional we attach to $\Delta_X$ on the left of (\ref{mainlemma})  contains too much curvature . The Mabuchi energy involves at most the Ricci curvature. In dimension $n\leq 2$ the only dually degenerate varieties are linearly embedded projective spaces (see \cite{tevelev}). Therefore the dual of a (nonlinear) projective curve is always a hypersurface, moreover the only curvature available is the Ricci curvature . Therefore the Main lemma applies, without modification, to space curves. In dimension two the second problem arises but not the first . Miraculously, Cayley's  $X$-\emph{hyperdiscriminant} (properly formatted)  eliminates both difficulties \emph{simultaneously}, moreover the hyperdiscriminant coincides with the usual discriminant in dimension one. 
Theorem A  follows from  working out the left hand side of (\ref{mainlemma}) in the case where $X$ \emph{has been replaced by} $\xhyp$. We always consider $\xhyp$ as a subvariety of $\mathbb{P}(M_{n\times (N+1)}^{\vee}(\mathbb{C}))$ via the Segre embedding. Observe that $G=\slnc$ acts on $M_{n\times (N+1)}^{\vee}(\mathbb{C})$  by the standard action on $\mathbb{C}^{N+1}$ and the trivial action on $\mathbb{C}^n$.

%%%%%%%%%%%%%%%%%%%%%%%%%%%%%%%%%%%%%%%%%%%%%%%%%%%%%%%%%%%%%%%%%%%%%%%%%%%%%%%%%%%%%%%%%%%%%%
\section{Completion of the proof of Theorem A }
 Theorem A follows at once from (\ref{mainlemma}) and the following proposition.
\begin{proposition}\label{mainprop}
 Let $X^n \hookrightarrow \cpn$ be a smooth, linearly normal algebraic variety of degree $d\geq 2$ . 
 Let   $R_{X}$  denote the \textbf{X-resultant} . Let $\Delta_{\xhyp}$ denote the \textbf{$X$-hyperdiscriminant}  .
 Then the Donaldson functional associated to the complex $K_{X\times \mathbb{P}^{n-1}}^{i}:= \bigwedge ^iJ_1(\mathcal{O}(1)|_{X\times \mathbb{P}^{n-1}})^{\vee}$ is given by
\begin{align} \label{bergm}
\begin{split}
& {\deg(R_X)}D_{\ K_{X\times \mathbb{P}^{n-1}}^{\bull} }(Ch_{n+1} ; H^{\bull}(\sigma),\ H^{\bull}(e) \ )=\nu_{\om}(\varphi_{\sigma})+ {\deg(\Delta_{\xhyp})}\log \frac {{||\sigma\cdot R_{X}||}^{2}}{||R_{X}||^2}  \ .
 \end{split}
\end{align}
\end{proposition}
 This entire section is devoted to the proof of this proposition. To begin let
\begin{align*}
w:=(w_1,w_2,\dots,w_n)\in \mathbb{C}^n \ra (1,T_1(w),T_2(w),\dots,T_N(w))\in \mathbb{C}^{N+1}
\end{align*}
be a local parametrization of $\tilde{X}$, where $\tilde{X}\subset \mathbb{C}^{N+1}\setminus\{0\} $ is the affine cone over $X$. Then $(w_1,w_2,\dots,w_n)$ are local coordinates on $X$. Observe that 
\begin{align*}
&e(w):=(1,T_1(w),T_2(w),\dots,T_N(w)) \\
\ \\
&f_i(z):=(0,\frac{\dl}{\dl w_i}T_1,\dots,\frac{\dl}{\dl w_i}T_N) \ .
\end{align*}
locally trivialize the bundle $J_1(\mathcal{O}(1)|_X)^{\vee}$. Note that $e(w)$ spans $\mathcal{O}(-1)|_X$ .  As we have remarked, the bundle of one jets is naturally a subbundle of the trivial bundle $X\times \mathbb{C}^{N+1}$. Let $\mathcal{I}_X=(F_{\alpha})$ denote a (finite) generating set for the homogeneous ideal of $X$. Then the jet bundle may be exhibited concretely as follows :
 \begin{align}\label{1jet}
J_1(\mathcal{O}(1)|_X)^{\vee} =\{(p,w)\in X\times \mathbb{C}^{N+1}|\ \nabla F_{\alpha}(p)\cdot w=0 \ \mbox{for all $\alpha$} \}\overset{\iota}{\hookrightarrow} X\times \mathbb{C}^{N+1}\ .
\end{align}
The (dual of) the jet bundle therefore inherits the standard Hermitian metric $h_{\mathbb{C}^{N+1}}$ from its embedding in $X\times \mathbb{C}^{N+1}$. 
We make extensive use of the following well known fact.
 \begin{proposition} There is an exact sequence of vector bundles on $X$.
 \begin{align}\label{xact}
 \begin{split}
 0\rightarrow \mathcal{O}_X(-1)\overset{\iota}{\rightarrow}&T^{1,0}(\tilde{X}) \overset{\pi}{\rightarrow} T^{1,0}(X)\otimes \mathcal{O}_X(-1)\rightarrow 0 
 \end{split}
 \end{align}
 \end{proposition}
 
 Since we will need an explicit description of the maps in what follows we recall the proof. Below we abuse notation as follows: on the one hand $\pi$ denotes the map
\begin{align*}
J_1(\mathcal{O}_X(1))^{\vee}\overset{\pi}{\rightarrow} T^{1,0}(X)\otimes \mathcal{O}_X(-1)\rightarrow 0 \ .
\end{align*}
On the other hand we \emph{also} denote by $\pi$ the projection onto $\cpn$
\begin{align*}
\pi:\mathbb{C}^{N+1}\setminus \{0\}\rightarrow \cpn \ .
\end{align*}
Finally we can define $\pi$ in (\ref{xact}) by the formula (where $\pi(v)=p$ )
\begin{align*}
J_1(\mathcal{O}_X(1))^{\vee}  \ni (p,w)\rightarrow \pi(p,w):= {\pi_{*}}|_v(w)\otimes v \in T^{1,0}(X)\otimes \mathcal{O}_X(-1)\ .
\end{align*}
The rationale for this follows from the fact that for all $w\in \mathbb{C}^{N+1}$ and $\alpha \in \mathbb{C}^*$ we have 
\begin{align*}
{\pi_{*}}|_{\alpha v}(w)=\frac{1}{\alpha}{\pi_{*}}|_{v}(w) \ .
\end{align*}

\begin{remark}
 \emph{On the one hand $\pi$ denotes the map
\begin{align*}
T^{1,0}(\tilde{X})\overset{\pi}{\rightarrow} T^{1,0}(X)\otimes \mathcal{O}_X(-1)\rightarrow 0 \ .
\end{align*}
On the other hand we \emph{also} denote by $\pi$ the projection onto $\cpn$
\begin{align*}
\pi:\mathbb{C}^{N+1}\setminus \{0\}\rightarrow \cpn \ .
\end{align*}
Finally we can define $\pi$ in (\ref{xact}) by the formula (where $\pi(v)=p$ )}
\begin{align*}
T^{1,0}(\tilde{X})  \ni (p,w)\rightarrow \pi(p,w):= {\pi_{*}}|_v(w)\otimes v \in T^{1,0}(X)\otimes \mathcal{O}_X(-1)\ .
\end{align*}
\end{remark}
Let $z\in\mathbb{C}^{N+1}\setminus \{0\} $  we define  
\begin{align}
g_{i\bar{j}}(z):=\frac{1}{|z|^4}(\delta_{ij}|z|^2-\bar{z}_iz_j ) \ . 
\end{align}

We define a Hermitian form $H_{i\bar{j}}(z)$ as follows:
\begin{align*}
H_{i\bar{j}}(z):=|z|^2g_{i\bar{j}}(z) \ .
\end{align*}
Then $H$ is a positive definite Hermitian form on $\mathcal{O}_X(-1)^{\perp}$. Moreover
\begin{align*}
h_{\mathbb{C}^{N+1}}|_{\mathcal{O}_X(-1)^{\perp}}=H \ .
\end{align*}
Therefore the standard Hermitian metric $h_{\mathbb{C}^{N+1}}$ descends to $\om\otimes h_{FS}$ on $T^{1,0}(X)(-1)$ .

Observe that $|e(w)|^2$  represents  the  Fubini Study local metric potential. Therefore the K\"ahler form on $X$ is given by
 
\begin{align*}
\om=\om_{FS}|_X=\frac{\sqrt{-1}}{2\pi}\dl_w\dlb_w\log|e(w)|^2 \ .
\end{align*}
Observe that
\begin{align*}
(e,f_i)=\frac{\dl}{\dl\bar{w}_i}|T|^2 
\quad
(f_i,f_j)=\frac{\dl^2}{\dl {w}_i\dl\bar{w}_j}|T|^2 \ .
\end{align*}
  By $f_i^{\perp}$ we mean the orthogonal projection of $f_i$ onto $\mathcal{O}(-1)^{\perp}$
\begin{align*}
f_i^{\perp}:=f_i-\frac{ (f_i\ ,\ e)}{|e|^2}e \ .
\end{align*}
Then with respect to the \emph{smooth} basis
\begin{align*}
\{e;\ f_1^{\perp}, f_2^{\perp}, \dots, f_n^{\perp}\}
\end{align*}
the matrix presentation of the metric $H$ has the shape
\begin{align*}
H_{\infty}=\begin{pmatrix}|e|^2&0&\dots&0\\
0&|e|^2g_{1\bar 1}&\dots&|e|^2g_{1\bar n}\\
0&|e|^2g_{2\bar 1}&\dots&|e|^2g_{2\bar n}\\
\dots&\dots&\dots&\dots\\
0&|e|^2g_{n\bar 1}&\dots&|e|^2g_{n\bar n}
\end{pmatrix}
\end{align*}
Let $H_{\mathcal{O}}$ be the matrix presentation of $H$ with respect to $\{e;\ f_1, f_2, \dots, f_n\}$, then it is easy to see that
$H_{\mathcal{O}}$ and $H_{\infty}$ are related by
\begin{align}\label{ortho}
H_{\mathcal{O}}=Q^TH_{\infty}\bar{Q} \ .
\end{align}
 The matrix $Q$ is given by
\begin{align}
\begin{split}
&Q=\begin{pmatrix}
1&q_1&q_2&\dots&q_n\\
\ \\
0&1&0&\dots&0 \\
\ \\
0&0&\dots&1&0\\
\ \\
0&0&\dots&0&1
\end{pmatrix}\\
\ \\
& q_i:=\frac{ (f_i\ ,\ e)}{|e|^2}\ .
\end{split}
\end{align}
Therefore 
\begin{align}
\det(H_{\infty})=\det(H_{\mathcal{O}})=|e|^{2(n+1)}\det(g_{i\bar j}(z)) \ .
\end{align}
 This gives the following pointwise identity of forms .    
\begin{align}\label{jet}
c_1(J_1(\mathcal{O}(1)|_X)^{\vee},\ h_{\mathbb{C}^{N+1}})= -(n+1)\om_{FS}|_X+\mbox{Ric}(\om_{FS}) \ .
\end{align}
(\ref{jet}) is a special case of a much more general ``metric splitting" of the Chern forms of the exact sequence
\begin{align} 
 \begin{split}
 0\rightarrow \mathcal{O}_X(-1)\overset{\iota}{\rightarrow}&  J_{1}(\mathcal{O}_X(1))^{\vee}\overset{\pi}{\rightarrow} T^{1,0}(X)\otimes \mathcal{O}_X(-1)\rightarrow 0 \ .
 \end{split}
 \end{align}
Since this is of such importance for the main result of this article and has played a significant role in the field in general we take time to discuss it.

Let $X$ be a complex manifold, and consider a short exact sequence of analytic vector bundles over $X$
\begin{align*}
0\ra \mathcal{E}_0\overset{j}{\ra} \mathcal{E}_1\overset{\pi}{\ra} \mathcal{E}_2\ra 0\ .
\end{align*} 
It is well known that the following identities are valid in $ H^{\bull}(X,\mathbb{C})$.
\begin{align}
\begin{split}
&c_{\tau}(\mathcal{E}_1)=c_{\tau}(\mathcal{E}_0)c_{\tau}( \mathcal{E}_2)     \\
\ \\
&Ch(\mathcal{E}_1)=Ch(\mathcal{E}_0)+Ch(\mathcal{E}_2) \ .
\end{split}
\end{align}
When the terms of the sequence are equipped with Hermitian metrics $h_j$ and corresponding curvatures $\Theta(\mathcal{E}_j\ , \ h_j)$ we may ask if the \emph{pointwise identities} hold
\begin{align}\label{pointwise}
\begin{split}
&\det\big(\tau I_{r_1}+\Theta(\mathcal{E}_1\ , \ h_1)\big)=\det\big(\tau I_{r_0}+\Theta(\mathcal{E}_0\ , \ h_0)\big)\det\big(\tau I_{r_2}+\Theta(\mathcal{E}_2\ , \ h_2)\big) \\
\ \\
&Tr\big(\exp(\Theta(\mathcal{E}_1\ , \ h_1))\big)=Tr\big(\exp(\Theta(\mathcal{E}_0\ , \ h_0))\big)+Tr\big(\exp(\Theta(\mathcal{E}_2\ , \ h_2))\big)\ .
\end{split}
\end{align}
In general they do not. An important example, that has in some sense shaped the field of K-stability, is the following. Let $X_F$ be a smooth hypersurface of degree $d\geq 2$ inside $\mathbb{P}^{n+1}$. Then we have the famous \emph{adjunction sequence}
 \begin{align}
 0\ra T^{1,0}_{X_F}\ra T^{1,0}_{\mathbb{P}^{n+1}}|_{X_F}\ra\mathcal{O}_{\mathbb{P}^{n+1}}(d)|_{X_F}\ra 0 \ .
 \end{align}
Equip each term $\mathcal{E}_{j}$ in the sequence with the induced Fubini-Study metric. In \cite{kenhyp} Tian has shown that
\begin{align}
Ric({\om_{FS}}|_{X_F})=(n+2-d){\om_{FS}}|_{X_F}-\frac{\sqrt{-1}}{2\pi}\dl\dlb\log\left(\frac{\sum_{0\leq j\leq n+2}|\frac{\dl F}{\dl z_j}|^2}{||z||^{2d-2}}\right) \ .
\end{align}
This shows that the pointwise identity already fails for $c_1$. This phenomena has been thoroughly analyzed by Bismut, Gillet, and Soul\'e in their famous 1988 paper ``\emph{Analytic Torsion and Holomorphic Determinant Bundles I}" (see \cite{bgs1}). The obstructions to splitting are called the \emph{secondary classes of Bott-Chern}. Precisely, Bismut, Gillet, and Soul\'e construct forms 
$\widetilde{Ch}(\mathcal{E}_{\bull}\ ; \ h_{\bull})$ which are unique modulo $\dl$ and $\dlb$ terms satisfying the following 
\begin{align}
\sum_{j=0}^2(-1)^jTr\big(\exp(\Theta(\mathcal{E}_j\ , \ h_j))\big)=-\frac{\sqrt{-1}}{2\pi}\dl\dlb \widetilde{Ch}(\mathcal{E}_{\bull}\ ; \ h_{\bull}) \ .
\end{align}
They construct similar classes for the total Chern class $c_{\tau}$. These secondary forms all have the property that, whenever the sequence splits as a \emph{holomorphic Hermitian} sequence the forms \emph{vanish identically}. 
 Since the jet complex does \emph{not} split metrically (where each term has its natural metric) a somewhat surprising fact about this complex is the following .
\begin{proposition} 
The Bott-Chern secondary classes of the jet complex with respect to the natural metrics vanish identically. Precisely, there is a pointwise identity of forms on $X$
 \begin{align}\label{jetpointwise}
\begin{split}
&\det\big(\tau I_{n+1}+\Theta(J_{1}(\mathcal{O}_X(1))^{\vee}  \ , \ h_{\mathbb{C}^{N+1}})\big)= \big(\tau  -\om_{FS}|_X\big)\det\big((\tau-\om_{FS}|_X)I_n + F_{\om_{FS}|_X} \big) \\
\ \\
&Tr\big(\exp(\Theta(J_{1}(\mathcal{O}_X(1))^{\vee}  \ , \ h_{\mathbb{C}^{N+1}}   )\big)= \exp(-\om_{FS}|_X)+Tr\big(\exp(     -\om_{FS}|_XI_n + F_{\om_{FS}|_X})\big)\ .
\end{split}
\end{align}
$F_{\om_{FS}|_X}$ denotes the full Riemann curvature tensor of $(X\ , \ \om_{FS})$ . 
\end{proposition}

A proof of (\ref{jetpointwise}) will be provided in the paragraph below for the \emph{top} Chern class $c_n$. This is all that is required
for my purpose.

Let
\begin{align*}
0\ra \mathcal{S}\overset{j}{\ra} \mathcal{E}\overset{\pi}{\ra} \mathcal{Q}\ra 0
\end{align*}
be a short exact sequence of holomorphic vector bundles on some  complex manifold $X$. Assume that $\mathcal{E}$ is equipped with a Hermitian metric $h$. Then $\mathcal{S}$ acquires a metric by restriction and $\mathcal{Q}$ by the (smooth) isomorphism
\begin{align*}
\mathcal{S}^{\perp}\overset{\pi}{\cong} \mathcal{Q} \ .
\end{align*}
The purpose of this paragraph is to analyze the curvature $F^{\mathcal{E}}$ in terms of $F^{\mathcal{S}}$ and $F^{\mathcal{Q}}$. This is due to Griffiths (see \cite{griffextII} section 6.3, \cite{griffherm} section 2 (d) ) .   This will then be applied to the jet sequence (\ref{xact}).   $D$ always denotes the unique holomorphic Hermitian connection.  
\begin{align}
\begin{split}
D^{\mathcal{E}}&=D^{\mathcal{E}}\circ j^*+D^{\mathcal{E}}\circ\pi^*\circ\pi\\
\ \\
&=D^{\mathcal{S}}\circ j^*+\pi^*\circ D^{\mathcal{Q}}\circ \pi+ \alpha\circ j^* +\beta\circ\pi^*\circ\pi \ .
\end{split}
\end{align}
Where we have defined
\begin{align}
\begin{split}
&\alpha:= D^{\mathcal{E}}\circ j^*-D^{\mathcal{S}}\circ j^* \\
\ \\
&\beta:=D^{\mathcal{E}}\circ \pi^*\circ\pi-\pi^*\circ D^{\mathcal{Q}}\circ \pi \ .
\end{split}
\end{align}
  $\alpha$ is the \emph{second fundamental form} of the inclusion $0\ra \mathcal{S}\overset{j}{\ra} \mathcal{E}$ . 
 
\begin{proposition}
\begin{align*}
&\alpha \in C^{\infty}(\Omega^{1,0}_X\otimes Hom(\mathcal{S},\ \mathcal{S}^{\perp}))\\
\\
&\beta=-\alpha^* \in C^{\infty}(\Omega^{0,1}_X\otimes Hom(\mathcal{S}^{\perp},\mathcal{S})) \ .
 \end{align*}
 \end{proposition}
In particular we have that
\begin{align}
\begin{split}
&\pi\circ \alpha \in C^{\infty}(\Omega^{1,0}_X\otimes Hom(\mathcal{S},\mathcal{Q}))  \\
\ \\
& \beta\circ \pi^*\in C^{\infty}(\Omega^{0,1}_X\otimes Hom(\mathcal{Q},\mathcal{S})) \\
\ \\
&(\pi\circ\alpha)\wedge(\beta\circ\pi^*)\in C^{\infty}(\Omega_X^{1,1}\otimes Hom(\mathcal{Q},\mathcal{Q})) \\
\ \\
&(\beta\circ\pi^*)\wedge(\pi\circ\alpha)\in C^{\infty}(\Omega_X^{1,1}\otimes Hom(\mathcal{S},\mathcal{S}))\ .
\end{split}
\end{align}
Then we have the basic \emph{curvature formula}.
 \begin{proposition}
 \begin{align}
\begin{split}
F^{\mathcal{E}}&= F^{\mathcal{S}}\circ j^*+\pi^*\circ F^{\mathcal{Q}}\circ \pi + \pi^*\dlb_{Hom(\mathcal{S},\mathcal{Q})}(\pi\circ \alpha)\circ j^*\\
\ \\
&+D^{1,0}_{Hom(\mathcal{Q},\mathcal{S})}(\beta\circ \pi^*)\circ \pi+
(\beta\wedge\alpha)\circ j \circ j^*+(\alpha\wedge\beta)\circ\pi^*\circ\pi \ .
\end{split}
\end{align}
\end{proposition}
Now we return to our situation, $\mathcal{S}=\mathcal{O}_X(-1)$, $\mathcal{E}=J_1(\mathcal{O}(1)|_X)^{\vee}$, and $\mathcal{Q}=T^{1,0}_X(-1)$.
In this case we have the identifications
\begin{align}
\begin{split}
&Hom(\mathcal{S},\mathcal{Q})\cong T^{1,0}_X \\
\ \\
&Hom(\mathcal{Q},\mathcal{S}) \cong \Omega^{1,0}_X \ .
\end{split}
\end{align}
The next proposition is crucial, it identifies the second fundamental form $\alpha$, in particular it shows that $\alpha$ is \emph{metric independent}. This proposition can be traced back to Atiyah (see \cite{atiyah58}), the author thanks an anonymous referee for pointing this out.
\begin{proposition}
\begin{align}
\pi\circ\alpha=dw_1\otimes\frac{\dl}{\dl w_1}+dw_2\otimes\frac{\dl}{\dl w_2}+\dots+dw_n\otimes\frac{\dl}{\dl w_n}\ .
\end{align}
Consequently,
\begin{align}
\beta\circ\pi^*=-\sum_{1\leq i,j\leq n}g_{i\bar j}(w)\ d\bar{w_j}\otimes dw_i \ .
\end{align}
Therefore $\alpha$ is holomorphic, $\beta$ is parallel and the curvature operator reduces to
\begin{align}
F^{\mathcal{E}}&= F^{\mathcal{S}}\circ j^*+\pi^*\circ F^{\mathcal{Q}}\circ \pi +  
(\beta\wedge\alpha)\circ j \circ j^*+(\alpha\wedge\beta)\circ\pi^*\circ\pi \ .
\end{align}
\end{proposition}
\begin{proof}
The proof is a straightforward computation. To begin
\begin{align}
D^{J_1(\mathcal{O}(1)|_X)^{\vee}}(e)=\om_{11}\otimes e+\sum_{2\leq j \leq n+1}\om_{j1}\otimes f_{j-1} \ .
\end{align}
The matrix of connection forms is given by the usual rule
\begin{align}
\om_{ij}=\sum_{1\leq k\leq n+1}h^{ki}\dl h_{jk} \ .
\end{align}
Therefore we have
\begin{align}
\begin{split}
&\om_{j1}=h^{kj}\frac{\dl h_{1k}}{\dl w_i}d w_i=h^{kj}h_{i+1 k}dw_i=\delta_{i+1 j}dw_i=dw_{j-1} \qquad (\mbox{$j\geq 2$}) \\
\ \\
&\om_{11}=0 \ .
\end{split}
\end{align}
By the same token
\begin{align}
\dl\log(|e|^2)=\frac{(f_1,\ e)}{|e|^2}dw_1+\frac{(f_2,\ e)}{|e|^2}dw_2+\dots+\frac{(f_n,\ e)}{|e|^2}dw_n \ .
\end{align}
Therefore
\begin{align}
\begin{split}
\alpha(e)&=(f_1-\frac{(f_1,\ e)}{|e|^2}e)\otimes dw_1+(f_2-\frac{(f_2,\ e)}{|e|^2}e)\otimes dw_2+\dots +(f_n-\frac{(f_n,\ e)}{|e|^2}e)\otimes dw_n\\
\ \\
&=f_1^{\perp}\otimes dw_1+f_2^{\perp}\otimes dw_2+\dots +f_n^{\perp}\otimes dw_n \ .
\end{split}
\end{align}
Since $\pi(f_j)=e\otimes \frac{\dl}{\dl w_j}$ we are done.
\end{proof}

%%%%%%%%%%%%%%%%%%%%%%%%%%%%%%%%%%ATIYAH%%%%%%%%%%%%%%%%%%%%%%%%%%%%%%%%%%%%%%%%%%%
From the above we have that
\begin{align}
\begin{split}
&(\beta\wedge\alpha)\circ j \circ j^*= \sum_{1\leq i,j\leq n}g_{i\bar j}(w)\ d{w_i}\wedge d\bar{w_j}  = {\om_{FS}}|_X\otimes I_{\mathcal{O}(-1)}\\
\ \\
&(\pi\circ\alpha)\wedge(\beta\circ\pi^*)=\begin{pmatrix}
-g_{1\bar{l}}(w)dw_1\wedge d\bar{w_l}&-g_{2\bar{l}}(w)dw_1\wedge d\bar{w_l}&\dots&-g_{n\bar{l}}(w)dw_1\wedge d\bar{w_l}\\
\ \\
-g_{1\bar{l}}(w)dw_2\wedge d\bar{w_l}&-g_{2\bar{l}}(w)dw_2\wedge d\bar{w_l}&\dots&-g_{n\bar{l}}(w)dw_2\wedge d\bar{w_l}\\
\ \\
\dots&\dots&\dots\\
-g_{1\bar{l}}(w)dw_n\wedge d\bar{w_l}&-g_{2\bar{l}}(w)dw_n\wedge d\bar{w_l}&\dots&-g_{n\bar{l}}(w)dw_n\wedge d\bar{w_l}
\end{pmatrix}
\end{split}
\end{align}
Where we sum over repeated indices.
Therefore,
\begin{align}
F^{J_1(\mathcal{O}(1)|_X)^{\vee}}=  \pi^*\circ \left(-{\om_{FS}}|_X\otimes I_{T^{1,0}}+F^{T^{1,0}_X}_{\om}\right)\circ \pi+(\alpha\wedge\beta)\circ\pi^*\circ\pi \ .
\end{align}
Since $\alpha$ takes values in $\mathcal{S}^{\perp}$ we have
\begin{align}
\pi^*\circ (\pi\circ \alpha)\wedge(\beta\circ\pi^*)\circ\pi=(\alpha\wedge\beta)\circ\pi^*\circ\pi \ .
\end{align}
At the center of a normal coordinate system the second fundamental form operator $S:=(\pi\circ \alpha)\wedge(\beta\circ\pi^*)$ takes the shape
\begin{align}S=
\begin{pmatrix}
-dw_1\wedge d\bar{w}_1&-dw_1\wedge d\bar{w}_2&\dots&-dw_1\wedge d\bar{w}_n\\
-dw_2\wedge d\bar{w}_1&-dw_2\wedge d\bar{w}_2&\dots&-dw_2\wedge d\bar{w}_n\\
\dots&\dots&\dots&\dots\\
-dw_n\wedge d\bar{w}_1&-dw_n\wedge d\bar{w}_2&\dots&-dw_n\wedge d\bar{w}_n\\
\end{pmatrix}
\end{align}
Observe that
\begin{align}
Tr((\pi\circ\alpha)\wedge(\beta\circ\pi^*))=-\om_{FS}|_X \ .
\end{align}
Therefore
\begin{align}\label{1jet2}
Tr(F^{J_1(\mathcal{O}(1)|_X)^{\vee}})=-(n+1)\om_{FS}|_X+\mbox{Ric}(\om_{FS}|_X) \ .
\end{align}
 (\ref{1jet2}) is consistent with (\ref{jet}). 
 \begin{lemma}
Let $F$ denote the full curvature tensor of $\om_{FS}|_X$.  
Then, for all $k\geq 1$ we have that
\begin{align}
\mbox{Trace}(F^kS)\equiv 0 \ .
\end{align}
\end{lemma}
The proof takes up one line.
Recall that at the center of a normal coordinate system $F$ has the shape
\begin{align}
F_{i\bar{j}}=-\dl\dlb g_{j\bar{i}} \ .
\end{align}
Therefore
\begin{align}
\begin{split}
&(-1)^{k+1}\mbox{Trace}(F^kS)=\\
\ \\
& \sum 
\frac{\dl^2 g_{i_1\bar{i}}}{\dl z_{p_1}\dl \bar{z}_{q_1}}\frac{\dl^2 g_{i_2\bar{i_1}}}{\dl z_{p_2}\dl \bar{z}_{q_2}}\dots\frac{\dl^2 g_{i_k\bar{i}_{k-1}}}{\dl z_{p_k}\dl \bar{z}_{q_k}}dz_{p_1}\wedge d\bar{z}_{q_1}\wedge  \dots \wedge dz_{p_k}\wedge d\bar{z}_{q_k}\wedge dz_{i_k}\wedge d\bar{z}_{i} \ .
\end{split}
\end{align}
That this sum is identically zero follows at once from the identities
\begin{align}
\begin{split}
\frac{\dl^2 g_{i_k\bar{i}_{k-1}}}{\dl z_{p_k}\dl \bar{z}_{q_k}}&=\frac{\dl^2 g_{p_k\bar{i}_{k-1}}}{\dl z_{i_k}\dl \bar{z}_{q_k}} \\
\ \\
dz_{p_k}\wedge d\bar{z}_{q_k}\wedge dz_{i_k}&=(-1)dz_{i_k}\wedge d\bar{z}_{q_k}\wedge dz_{p_k} \ .
\end{split}
\end{align}
The proof is complete. $\Box$

The definition of $S$ implies at once that 
\begin{align}
S^2=\om S \ .
\end{align}
Since $\mbox{Trace}(S)=-\om$ we have the following
\begin{align}\label{pwrsum}
\mbox{Trace}(F+S)^k=\mbox{Trace}(F^k)-\om^k \ .
\end{align}
\begin{lemma}\label{sigma}
For any $A\in M_n(\mathbb{C})$ let $\sigma_k(A)$ denote the $kth$ elementary symmetric function of $A$. Then
\begin{align}
\sigma_k(F+S)=\sum_{j=0}^k(-1)^j\sigma_{k-j}(F)\om^j \ .
\end{align}
\end{lemma}
\begin{proof}
The identity obviously holds when $k=1$. We proceed by induction. So assume the identity for $1\leq j\leq k-1$. Newton's formula relating $\sigma_{k}(A)$ and $p_k(A):=\mbox{Trace}(A^k)$ together with (\ref{pwrsum}) imply that
\begin{align}
\begin{split}
&\sigma_{k}(F+S)=\frac{1}{k}\big\{\sum_{j=1}^k(-1)^{j+1}\sigma_{k-j}(F+S)(p_{j}(F)-\om^j)\big\} \\
\ \\
&\frac{1}{k}\big\{\sum_{j=1}^k\sum_{i=0}^{k-j}(-1)^{i+j+1}\sigma_{k-j-i}(F)p_j(F)\om^i\big\} 
+\frac{1}{k}\big\{\sum_{j=1}^k\sum_{i=0}^{k-j}(-1)^{i+j}\sigma_{k-j-i}(F)\om^{i+j}\big\} \ .
\end{split}
\end{align}
 
Rearrangement shows that
\begin{align}
 \begin{split}
& \sum_{j=1}^k\sum_{i=0}^{k-j}(-1)^{i+j+1}\sigma_{k-j-i}(F)p_j(F)\om^i  =\sum_{i=0}^{k-1}(-1)^i(k-i)\sigma_{k-i}(F)\om^i \\
\ \\
& \sum_{j=1}^k\sum_{i=0}^{k-j}(-1)^{i+j}\sigma_{k-j-i}(F)\om^{i+j}  =\sum_{i=1}^{k-1}(-1)^ii\sigma_{k-i}(F)\om^i \ .
\end{split}
\end{align}
Adding these two completes the proof of the proposition.
\end{proof}

\begin{corollary}\label{pointwise}There is a \textbf{pointwise} identity of differential forms
\begin{align}
c_n(J(\mathcal{O}_X(1))^{\vee}\ ,\ \om)=c_n(T^{1,0}_X(-1)\ , \ \om)-c_{n-1}(T^{1,0}_X(-1) \ , \ \om)  \om \ .
 \end{align}
\end{corollary}
\begin{proof}
To begin we have that
\begin{align}
\begin{split}
&c_n(T^{1,0}_X(-1)\ , \ \om)-c_{n-1}(T^{1,0}_X(-1) \ , \ \om)\om=\\
\ \\
&\sum_{j=0}^n(-1)^jc_{n-j}(F)\om^j+\sum_{j=0}^{n-1}(-1)^{j+1}(j+1)c_{n-j-1}(F)\om^{j+1}= \\
\ \\
&\sum_{j=0}^n(-1)^j(j+1)c_{n-j}(F)\om^j \ .
\end{split}
\end{align}
 By definition,
\begin{align}
\begin{split}
c_n((J(\mathcal{O}_X(1))^{\vee}\ ,\ \om)&=\sigma_n(\pi^*\circ  \{-{\om_{FS}}|_X  I_{T^{1,0}}+F^{T^{1,0}_X}_{\om}  + S\}\circ\pi ) \\
\ \\
&=\det \big(-{\om_{FS}}|_X  I_{T^{1,0}}+F^{T^{1,0}_X}_{\om}  + S\big)\ .
\end{split}
\end{align}
By lemma (\ref{sigma}) we have
\begin{align}
\sum_{j=0}^n(-1)^j\om^j\sigma_{n-j}(F+S)=\sum_{k=0}^n\sum_{j=0}^k(-1)^{n-(k-j)}\om^{n-(k-j)}c_{k-j}(F) \ .
\end{align}
Now the corollary amounts to the following
\begin{claim}
\begin{align}
\sum_{k=0}^n\sum_{j=0}^k(-1)^{n-(k-j)}\om^{n-(k-j)}c_{k-j}(F)=\sum_{i=0}^n(-1)^i(i+1)c_{n-i}(F)\om^i \ .
\end{align}
\end{claim}
The proof of the claim is similar to the proof of lemma (\ref{sigma}) and is left to the reader. This completes the proof of  corollary (\ref{pointwise}) .
 \end{proof} 
  
  Let $\xi\in \mathfrak{sl}(N+1,\mathbb{C})$ and let $\sigma=\exp(\xi)\in \slnc $. We introduce a one parameter family of metrics $H_t=(.\ ,\ .)_t$ on $J_{1}(\mathcal{O}_X(1))^{\vee}$ 
joining $h_{\mathbb{C}^{N+1}}=H_0$ to $H_{\sigma}=H_1$
by the rule
\begin{align}
\begin{split}
(V,W)_t:=(\exp(t\xi)W,\ \exp(t\xi)V) \quad V,W\in \mathbb{C}^{N+1}\  .
\end{split}
\end{align}
Then
\begin{align}
\begin{split}
&H_t|_{\mathcal{O}(-1)}=\exp(\vpt)|\ \cdot |^2 \\
\ \\
&H_t|_{\mathcal{O}(-1)^{\perp}}=\exp(\vpt)|\ \cdot |^2\otimes \om_{t} \ .
\end{split}
\end{align} 
Where $\vpt$ and $\om_t$ are given by
\begin{align}
\begin{split}
&\vpt:= \log\left(\frac{|\exp(t\xi)T|^2}{|T|^2}\right) \\
\ \\
&\om_t:= \om_{FS}|_X+\frac{\sqrt{-1}}{2\pi}\dl\dlb\vpt \ .
\end{split}
\end{align}

Our aim is to compute, for a \emph{general} $X$, the Donaldson functional
\begin{align}
D _{K^{\bull}}(Ch_{n+1}\ ; H^{\bull}(\sigma),\ H^{\bull}(e) \ ) 
\end{align}
with respect to the path $h_t$ . Then we will replace $X$ with $\xhyp$.

Another application of the Borel-Serre identity gives
\begin{align}\label{det}
 \begin{split}
& D _{K^{\bull}}(Ch_{n+1}\ ; H^{\bull}(\sigma),\ H^{\bull}(e) \ )=\\
\ \\
&(-1)\int_0^1\int_X\frac{\dl}{\dl b}\det\big(\pi^*\circ\{F_{\om_t|_X} -\om_t|_XI_n+S(t)\}\circ \pi+bU(t)\big)|_{b=0}dt \ .
\end{split}
\end{align}
Where $F_{\om_t|_X}$ is the full Riemann curvature tensor of $\om_{t}$ and $U(t)$ is the endomorphism
\begin{align}
U(t):=\big(\frac{d}{ds}H_s\cdot {H_s}^{-1} \big)^T|_t \ .
\end{align}
Computation of the determinant in (\ref{det}) at the point $(o\ , \ t)$ with respect to the local analytic frame 
\begin{align*}
 \{e=T\ , \  f_i-\frac{ (f_i\ ,\ e)_t}{|e|_t^2}|_o\ e  \} \ ,
\end{align*}
shows at once that
\begin{align}
\begin{split}
\frac{\dl}{\dl b}\det\big(\pi^*\circ\{F_{\om_t|_X}-\om_t|_XI_n+S(t)\}\circ \pi+bU(t)\big)|_{b=0}&=\dot{\vpt}\det\big(F_{\om_t|_X}-\om_t|_XI_n+S(t)\big) \\
\ \\
&=\dot{\vpt}c_n(J_1(\mathcal{O}(1)|_{X})^{\vee}\ ; \ h_t) \ .
\end{split}
\end{align}

The next proposition seems to have been known to Cayley, a modern proof has been provided by Weyman and Zelevinsky. We give a new proof of the result of Weyman and Zelevinsky based on the theorem of Beltrametti, Fania, and Sommese mentioned in section 2 . The ingredients of the proof are required at a later stage in our argument.
\begin{proposition}
Let $X\hookrightarrow \mathbb{P}^{N}$ be a smooth linearly normal subvariety of degree $d\geq 2$. 
 Let $\mu$ denote the average of the scalar curvature.
Then the hyperdiscriminant of format $(n-1)$ is well formed with degree given by
\begin{align}
\deg(\Delta_{\xhyp})=n(n+1)d-d\mu \ .
\end{align}
   In particular,    
 \begin{align*}
\begin{split}
&\deg\Big(\Delta_{X^n_d\times\mathbb{P}^{n-1}}\Big)=n(n+1)(d-1) \qquad \mbox{( $X^n_d$ is the $dth$ Veronese image on $\mathbb{P}^n$ . )} \\
\ \\
&\deg\Big(\Delta_{X\times \mathbb{P}^{n-1}}\Big)=n\prod_{i=1}^kd_i\Big(\sum_{i=1}^kd_i-k\Big) \qquad \mbox{( $X\subset \mathbb{P}^{n+k+1}$ a complete intersection . )}\\
\ \\
&\deg\Big(\Delta_X\Big)=2d-2+2g \qquad \mbox{( $X$ a smooth curve of genus $g$ . )}
\end{split}
\end{align*}
\end{proposition}
\begin{proof}
Recall the \emph{smooth} isomorphism 
\begin{align}\label{splitting}
\Omega^{1,0}_{\mathbb{P}^{n-1}}\oplus \mathcal{O}\cong \overbrace{\bigoplus \mathcal{O}(-1)}^{n}\ .
\end{align}
The short exact sequence
\begin{align}
0\ra \Omega^{1,0}_{X\times \mathbb{P}^{n-1}}(1)\ra J_1(\mathcal{O}(1)|_{X\times \mathbb{P}^{n-1}})\ra\mathcal{O}(1)|_{X\times \mathbb{P}^{n-1}}\ra 0 \ 
\end{align}
 implies the Chern class identity
\begin{align}
c( J_1(\mathcal{O}(1)|_{X\times \mathbb{P}^{n-1}}))=c(\Omega^{1,0}_{X\times \mathbb{P}^{n-1}}(1))(1+\om_{FS}+\om) \ .
\end{align}
Recall that the restriction of the hyperplane from the Segre embedding  of $X\times \mathbb{P}^{n-1}$ is the tensor product   
\begin{align}
\mathcal{O}_{\mathbb{P}^{(N+1)n-1}}(1)|_{X\times \mathbb{P}^{n-1}}\cong \mathcal{O}(1)|_{X}\otimes  \mathcal{O}_{\mathbb{P}^{n-1}}(1) \ .
\end{align}
Next we have the obvious holomorphic splitting 
\begin{align}
\Omega^{1,0}_{X\times \mathbb{P}^{n-1}}(1)\cong \Omega^{1,0}_{X}(1)\otimes \mathcal{O}_{\mathbb{P}^{n-1}}(1)\oplus \Omega^{1,0}_{\mathbb{P}^{n-1}}(1)\otimes
\mathcal{O}(1)|_{X} \ .
\end{align}
Therefore
\begin{align}
c(\Omega^{1,0}_{X\times \mathbb{P}^{n-1}}(1))=c(\Omega^{1,0}_{X}(1)\otimes \mathcal{O}_{\mathbb{P}^{n-1}}(1))c(\Omega^{1,0}_{\mathbb{P}^{n-1}}(1)\otimes
\mathcal{O}(1)|_{X}) \ .
\end{align}
 By (\ref{splitting}) we have the smooth isomorphism over $X\times \mathbb{P}^{n-1}$ 
\begin{align}
\Omega^{1,0}_{\mathbb{P}^{n-1}}(1)\otimes\mathcal{O}(1)|_{X}\oplus  \mathcal{O}_{\mathbb{P}^{n-1}}(1)\otimes \mathcal{O}(1)|_{X}\cong\overbrace{\bigoplus \mathcal{O}(1)|_{X}}^n
\end{align}
Taking the total Chern class then gives
\begin{align}
c(\Omega^{1,0}_{\mathbb{P}^{n-1}}(1)\otimes
\mathcal{O}(1)|_{X})(1+\om_{FS}+\om)=(1+\om_{FS})^n \ .
\end{align}
Therefore we have that
\begin{align}
c( J_1(\mathcal{O}(1)|_{X\times \mathbb{P}^{n-1}}))= c(\Omega^{1,0}_{X}(1)\otimes \mathcal{O}_{\mathbb{P}^{n-1}}(1))(1+\om_{FS})^n\ .
\end{align}
Next we require the well known identity. Let $E$ be a rank $r$ vector bundle and $L$ a line bundle and $0\leq p\leq r$ an integer, then
\begin{align}
c_p(E\otimes L)=\sum^p_{i=0}\binom{r-i}{p-i}c_i(E)c_1(L)^{p-i} \ .
\end{align}
We see that
\begin{align}
\begin{split}
c_{n-1}(\Omega^{1,0}_{X}(1)\otimes \mathcal{O}_{\mathbb{P}^{n-1}}(1))&=\binom{n}{n-1}\om^{n-1}+O(\om^{n-2}) \\
\ \\
c_n(\Omega^{1,0}_{X}(1)\otimes \mathcal{O}_{\mathbb{P}^{n-1}}(1))&=c_1(\Omega^{1,0}_{X}(1))\om^{n-1}+O(\om^{n-2}) \\
\ \\
&=c_1(\Omega^{1,0}_{X})\om^{n-1}+n\om_{FS}\om^{n-1}+O(\om^{n-2}) \ .
\end{split}
\end{align}
So we see that 
\begin{align}
c( J_1(\mathcal{O}(1)|_{X\times \mathbb{P}^{n-1}}))= \big(c_1(\Omega^{1,0}_{X})\om^{n-1}+n\om_{FS}\om^{n-1}+n\om^{n-1}+O(\om^{n-2})\big)(1+\om_{FS})^n
\end{align}
From this the component of top dimension is easily seen to be
\begin{align}\label{topform}
\begin{split}
c_{2n-1}( J_1(\mathcal{O}(1)|_{X\times \mathbb{P}^{n-1}}))&= nc_1(\Omega^{1,0}_{X})\om^{n-1}_{FS} \om^{n-1}+n^2\om^n_{FS}\om^{n-1}+n\om_{FS}^n\om^{n-1}\\
\ \\
&= nc_1(\Omega^{1,0}_{X})\om^{n-1}_{FS} \om^{n-1}+n(n+1)\om^n_{FS}\om^{n-1}\ .
\end{split}
\end{align}

Next we show that the integral
\begin{align}
\int_X c_{2n-1}( J_1(\mathcal{O}(1)|_{X\times \mathbb{P}^{n-1}}))=n(n+1)d-d\mu >0
\end{align}
if and only if $d\geq 2$. The proof is a simple excercise in the adjunction formula, shown to the author by Lev Borisov.
To begin, let $H_1,H_2,\dots,H_{n-1}$ be  generic hyperplanes in $\cpn$. Let $\mathcal{C}_g$ denote the intersection
\begin{align}
\mathcal{C}_g:= \cap_{1\leq j\leq n-1}H_j\cap X \ .
\end{align}
Then $\mathcal{C}_g$ is a smooth curve of genus $g$. Let $K$ denote the canonical bundle of $\mathcal{C}_g$, then
\begin{align}
 2g-2=\int_{\mathcal{C}_g}c_1(K) \ .
 \end{align}
There is an exact sequence
\begin{align}
0\ra T^{1,0}_{\mathcal{C}_g}\ra T^{1,0}_{X}|_{\mathcal{C}_g}\ra\overbrace{\bigoplus \mathcal{O}_{\cpn}(1)|_{\mathcal{C}_g}}^{n-1}\ra 0
\end{align}
from which we deduce the isomorphism 
\begin{align}
 \mathcal{O}_{\cpn}(n-1)\otimes K_{X} \cong K \ .
\end{align}
Therefore,
\begin{align}
\begin{split}
2g-2&=\int_X\big(-\mbox{Ric}(\om |_X)+(n-1)\om \big)\om^{n-1} \\
\ \\
&=-\frac{d\mu}{n}+d(n-1) \ .
\end{split}
\end{align}
Since $g\geq 0$ and $d\geq 2$ we have the inequalities 
\begin{align}
0\leq n(n-1)d+2n-d\mu \leq n(n-1)d+dn-d\mu <n(n+1)d-d\mu \ .
\end{align}
 That the hyperdiscriminant is well formed, as well as the degree formula follow at once. 
 \end{proof}
From our previous work on the pointwise splitting of the Chern forms \footnote{We have made tacit use of the fact that the splitting holds for the Fubini-Study metric on $\mathbb{P}^n$. Precisely $c(T^{1,0}_{\mathbb{P}^n} \ ; \ \om)=(1+\om)^{n+1}$ \emph{pointwise}.}  and (\ref{topform}) we have the following 
\begin{claim}
There is a \textbf{pointwise identity} of forms on $\xhyp$
\begin{align}
c_{2n-1}( J_1(\mathcal{O}(1)|_{X\times \mathbb{P}^{n-1}})\ ; \ H_t)= -n\mbox{Ric}(\om_{\vpt}){\om^{n-1}_{\vpt}}  \om^{n-1}+n(n+1)\om_{\vpt}^n \om^{n-1}\ .
\end{align}
\end{claim}
We sum up the result of our work in the following
\begin{proposition}
Let $K^{\bull}=K_{X\times\mathbb{P}^{n-1}}^{\bull} $ denote the Cayley-Koszul associated to $\xhyp$ in its Segre embedding. Then
\begin{align}
\begin{split}
&D _{K^{\bull}}(Ch_{n+1}\ ; H^{\bull}(\sigma),\ H^{\bull}(e) \ )=\\
\ \\
&\int_0^1\int_{\xhyp}\dot{\vpt}\big\{-n\mbox{Ric}(\om_{\vpt})\om^{n-1}_{\vpt}  \om^{n-1}+n(n+1)\om_{\vpt}^n \om^{n-1}\big\} \ .
\end{split}
\end{align}
\end{proposition}
Next we require the following well known result which follows easily from the techniques used in this paper.\newline
\begin{theorem}\label{chownorm} [Tian \cite{kenhyp}, Zhang \cite{zhang}, Paul \cite{gacms}]   
  Let $X$ be an n-dimensional subvariety of \cpn, and let $R_{X}$ denote the $X$-resultant.
Then
\begin{align}\label{chow}
\begin{split}
&-{\deg}(X)(n+1)F_{\omega}^{0}(\varphi_{\sigma})= \log  \frac{{||\sigma\cdot R_{X} ||}^{2}}{{||R_{X}||}^{2}}  \ ; \quad
B:= \mathbb{P}H^{0}(\mathbb{G},\mathcal{O}(d)) \ .
\end{split}
\end{align}  
  $\mathbb{G}:= \mathbb{G}(N-n-1,\mathbb{P}^{N})$ denotes the Grassmannian of $N-n-1$ linear subspaces of $\cpn$, and the energy $F_{\omega}^{0}(\varphi)$ is defined as follows
  \begin{align}
F_{\omega}^{0}(\varphi):= -\int_0^1\int_X \dot{\vpt}\frac{\om^n_{\vpt} }{V} \ .
\end{align}

\end{theorem}
 
This ends the proof of proposition (\ref{mainprop}). The proof of the Theorem A is now complete.  
Theorem B follows at once from Theorem A and proposition (\ref{enasymp}) , Theorem C follows from proposition (\ref{kempfnessP}), we formulate an apparently stronger but equivalent form of Theorem D as follows. \\
\ \\
\noindent{\textbf{Theorem D (strong form) .}} \emph{Let $X\ra \cpn$ be a smooth, linearly normal algebraic variety of degree $d\geq 2$. Then $X$ is K-stable if and only if for all maximal algebraic tori $H$ and all $m>>0$ there is a constant $C=C(H)>0$ such that}
\begin{align}
\nu_{\om}(\vp_{\tau})\geq \frac{\deg(\Delta)\deg(R)}{m} J_{\om}(\vp_{\tau})-C \quad \tau \in H \ .
\end{align}
  \\
\ \\

The strong form of Theorem D follows at once from Theorem A, proposition (\ref{kempfnessproper}), Theorem (\ref{chownorm}) and Sun's lemma (proposition (\ref{sun'slemma}) ) .  Theorem E also follows from proposition (\ref{sun'slemma}) and  Theorem F is a consequence of corollary (\ref{generic}) . Corollary (\ref{ke}) part i) follows from Theorem C and \cite{chentian08} and part iii) requires Tian's properness theorem (\ref{tian97}) . Corollary (\ref{examples}) parts i) and ii) follow from proposition (\ref{kempfnessproper}), Theorem A, and the remark immediately following definition (\ref{proper}) .  We leave further details to the reader.  $\Box$
%%%%%%%%%%%%%%%%%%%%%%%%%%%%%%%%%%%%%%%%%%%%%%%%%%%%%%%%%%%%%%%%%%%%%%%%%%%
 \begin{center}{\textbf{Acknowledgements}}\end{center}
   The ideas of Gang Tian have been indispensible, the whole architecture of my program was inspired by his work on special degenerations and the generalized Futaki invariant. In a sense this article is the result of carefully studying his 1994 paper ``\emph{The K-Energy on Hypersurfaces and Stability}" \cite{kenhyp}. I thank him for his willingness to share his ideas with me over the many years that we have known one another and for being so supportive during the difficult time I faced during my postdoctoral studies. The author was inspired through separate conversations with Eckart Viehweg, who suggested projective duality in the summer of 2007 while the author visited with him in Essen, and Bernd Sturmfels, who, during a visit to Madison, emphasized that I should seek out the relevant polytopes. I thank Xiuxiong Chen, Jeff Viaclovsky, Sun Song, Lev Borisov, Yanir Rubinstein, Jordan Ellenberg and Chi Li for many helpful conversations. Sun in particular suggested proposition (\ref{sun'slemma}) which led to the weak numerical criterion for the Mabuchi energy, and a remark of Ellenberg's led to Theorem F . I thank my colleagues in the mathematics department at Madison for providing an outstanding environment for research. We should also mention that the excellent survey by Tevelev \cite{tevelev} has been very useful and instructive. Finally, the author owes an enormous debt to Arthur Cayley whose spectacular contributions he was made aware of through the fundamental work of Gelfand, Kapranov and Zelevinsky \cite{gkz}.
\bibliography{ref.bib}
\end{document}